%% file: article.tex
\documentclass{article} 
\usepackage{iclr2020_conference,times}

\usepackage{url}
\usepackage{amsmath}
\usepackage{amssymb}
\usepackage{amsthm}
\usepackage{mathtools}
\usepackage{float}
\usepackage{tikz}
\usetikzlibrary{arrows,snakes,backgrounds}
\usepackage{graphicx}
\usepackage{caption}
\usepackage{subcaption}
\usepackage{pgfmath}
\usepackage[boxed,resetcount]{algorithm2e}

 \title{FMM-Net: neural network architecture based on the fast multipole method}

\iclrfinalcopy
\author{
    Daria Sushnikova \\
    Skolkovo Institute of Science and Technology \\
    Moscow, Russia, 143026 \\
    \texttt{d.sushnikova@skoltech.ru} \\
\And
    Pavel Kharyuk \\
    Skolkovo Institute of Science and Technology \\
    Moscow, Russia, 143026 \\
    Marchuk Institute of Numerical Mathematics of the Russian Academy of Sciences \\
    Moscow, Russia, 119991 \\
    \texttt{kharyuk.pavel@gmail.com} \\
\And
    Ivan Oseledets \\
    Skolkovo Institute of Science and Technology \\
    Moscow, Russia, 143026 \\
    Marchuk Institute of Numerical Mathematics of the Russian Academy of Sciences \\
    Moscow, Russia, 119991 \\
    \texttt{ivan.oseledets@gmail.com} \\
}

\newtheorem{definition}{Definition}[section]
\newtheorem{example}{Example}[section]
\newtheorem{rmk}{Remark}

\newcommand{\size}{0.7}
\input{figs/close_mat.tex} 

\begin{document}

\maketitle

\begin{abstract}
    In this paper, we propose a new neural network architecture
    based on $\mathcal{H}^2$~matrix. Even though networks with
    $\mathcal{H}^2$-inspired architecture already exists, and our approach is designed
    to reduce memory costs and improve performance by taking into account the
    sparsity template of $\mathcal{H}^2$~matrix. In numerical comparison with
    alternative neural networks, including the known $\mathcal{H}^2$~based
    ones, our architecture showed itself as beneficial in terms of performance,
    memory, and scalability.
\end{abstract}

\section{Introduction}
\label{sec-1}

Neural networks (NN) have become extremely popular in the recent decade.
One of the magistral directions in this area is a search for the new architectures
that suit the specific problems. To name a few: UNet was proposed for
segmentation of medical images \citep{ronneberger-unet-2015}; deep convolutional
neural networks such as AlexNet \citep{krizhevsky-alexnet-2012}, VGGNet
\citep{simonyan-vgg-2014}, ResNet \citep{he-resnet-2015} designed to work with
visual data; bidirectional recurrent neural networks for the sequence to sequence
tasks \citep{schuster1997bidirectional,graves2013speech}.

There are many problems where it is required to approximate mapping
from one function to another, and it is necessary to design problem-specific architectures of neural networks that are adapted to the
structure of the problem. Such structures often arise from partial
differential equations (PDE) and integral equations (IE). If one
has an access to the
sequence of right-hand sides and corresponding solutions, the integral operator may
be learned on this basis. The only question is to define the structure of the
model which approximates the integral operator. It is well-known that integral
operators can be approximated well by block low-rank matrices with hierarchical
structure~\citep{GrRo-fmm-1987,hackbusch-h-1999,tee-mosaic-1996}, and one may
impose a such constraint on the learning model. Hence we came to the idea to shape
the neural network architecture as an extension of hierarchical matrix structure,
particularly, the $\mathcal{H}^2$-matrix ~\citep{hackbusch-h2-2000,Borm-h2-2010}.

In this approach, the forward pass of the neural network is similar to
$\mathcal{H}^2$~matrix by vector multiplication. The proposed architecture inherits
the double-tree structure of $\mathcal{H}^2$~matrices and may be considered as
a branching network, like Inception network~\citep{szegedy-inception-2017}.
However, the summation performed according to the second tree makes the
architecture similar to ResNet~\citep{he-resnet-2015}. The resulting network
naturally shares the scaling property of $\mathcal{H}^2$~matrices that may
contribute to its performance.

It is worth noting that $\mathcal{H}^2$~matrix structure is dependent on how
one selects close and far regions. Compared to other work on this theme
\citep{fan-hmix-2018,fan-h2mix-2019} where authors used band matrix for close
and tree-to-tree transfer matrices, we considered them as the block-sparse
matrices. Our architecture extensively relies on the sparsity template
resulting in replacing separate blocks of block-sparse matrices with the sequence of feed-forward fully-connected layers of small size.

In Section~\ref{sec-arc}, we briefly discuss the structure of the
$\mathcal{H}^2$~matrix and describe the architecture
of the neural network based on the $\mathcal{H}^2$~matrix and details of its
implementation. Section~\ref{sec:exp} contains experiments with implemented
neural network and a comparison of the proposed neural network with the reference
networks. Conducted experiments demonstrate the practical benefits of the proposed
architecture and provide a comparison with other existing approaches.

\section{Neural network architecture based on $\mathcal{H}^2$~matrix}
\label{sec-arc}

\subsection{From Fast Multipole Method (FMM) to FMMNet}

Let us consider two sets of points (point clouds), $x = \{x_i\}_{i = 1 
\dots N}$ and $y = \{ y_i\}_{j = 1 \dots N}$, $x_i, y_j \in
\mathbb{R}^{d}$, where $d=1, 2, 3$. We are referring to them
as source and receiver points respectively. Let us assume that
two functions, $q(x)$ and $w(y)$ are defined in subspaces
$X, Y \subset \mathbb{R}^{d}$ so that $x_i \in X$, $y_j \in
Y$, and $q_i = q(x_i)$, $w_j = w(x_j)$ for short. Consider
the mapping $\phi: q(X) \to w(Y)$ with $q(X) = \{h_x | h_x =
q(x) \, \forall x \in X\}$ and $w(Y) = \{h_y | h_y =
w(y)\, \forall y \in Y\}$, and our task is to approximate it.
The simplest approximation in this regard is a linear
transformation $w_j = \sum_{i=1}^{N} a_{ij} q_i$ with weights $a_{ij}$ which are interactions between $i$-th and $j$-th particles; $q(x)$, $w(y)$ and weights $a_{ij}$ are
dependent on the problem. One example is the N-body problem,
where $q_i$ is a mass of point $x_i$, $w_j$ is a value of
force in point $y_j$, and $a_{ij}$ is an inverse
distance between points $x_i$ and $y_j$, $a_{ij} = 1/\rho(x_i,
y_j)$. Another example is a discretization of
some integral equation with a smooth kernel.

Computation of such sum takes $O(N^2)$ operations. However,
due to the special structure of the problem, we can approximately
evaluate it using the Fast Multipole Method (FMM) in $O\left(N
\ln\left( 1/\varepsilon \right) \right)$ within a tolerance
$\varepsilon$. Detailed information can be found in \cite{}, 
but here we are interested only in the computational graph of the
algorithm to use it as a basis for neural network 
architecture.

\begin{figure}[H]
\centering
\input{figs/tree.tex}
\caption{Block cluster trees within FMM matrix by vector product. Each node of a tree denotes vectors, and each junction of nodes corresponds to a linear transformation of them.}
\label{fig:tree}
 \end{figure}
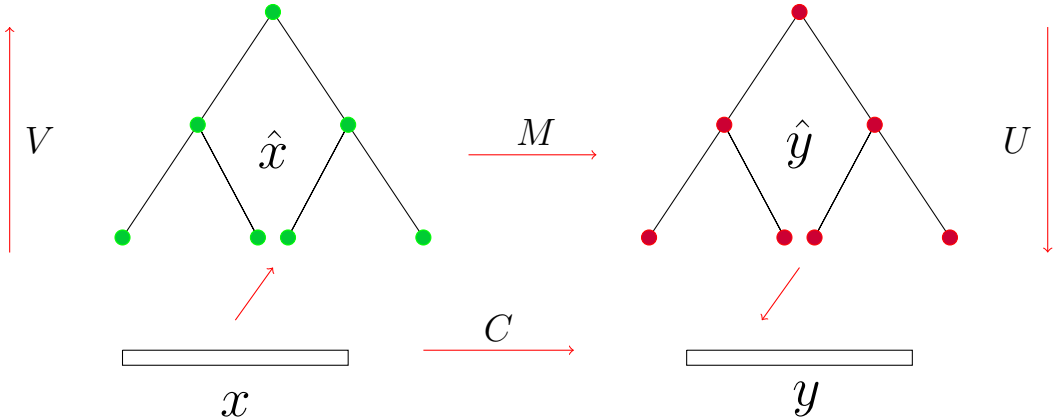
 
The construction of such architecture consists of 2 steps:
construction of two block cluster trees according to FMM and
generalization of intermediate linear mappings. At first we process
geometric information in order to create two cluster trees,
and then basing on so-called
separability criteria we create a block cluster tree which is
stored as a pairwise interaction list. There are many ways to do
that, see illustration ... [Huckbush, Tyrtyshnikov, Roklin].

To make the second step, we turn our attention to algebraic
interpretation of FMM. It can be written using the formalism of  $\mathcal{H}^2$-matrices in the following form:

\begin{equation}
    \label{eq:h2}
    w = Cq + \sum_{l=L}^{1} \left( U_L \ldots U_{L-l+1} \right) M_l
    \left( V_{L-l+1} \ldots V_{L} \right) q,
\end{equation}

where $C \in \mathbb{R}^{N \times N}$ is referred to as close
interaction matrix, $U_l, V_l$ are interpolator from and projector
on the rough grid (transition between so-called levels of
decomposition), $M_l$ is a mapping from $q(X)$ to $w(Y)$ on the
level $l$, $l=\overline{1,L}$, $w$ and $q$ are defined on two cloud
points $X$ and $Y$. All matrices included into
decomposition~\ref{eq:h2} have specific structure: $U_l$ and $V_l$
are block-diagonal, $C$ and $M_l$ are block-sparse. The sparsity
the pattern of the latter ones is defined by the structure of a certain
problem. The corresponding computational graph is presented in
Figure~\ref{fig:hnn}.

\begin{figure}[H]
\centering
    \resizebox{0.6\textwidth}{!}{
        \input{figs/nn.tex}
    }
\caption{Computational diagram of FMM matrix by vector product.
Mappings $G_0$, $G_l$, $l=\overline{1, L}$ are linear transformations with matrices $C$, $M_l$ respectively. In case
of FMMNet these mappings are replaced by nonlinear mappings
parametrized by neural networks with a certain structure.}
\label{fig:hnn}
\end{figure}
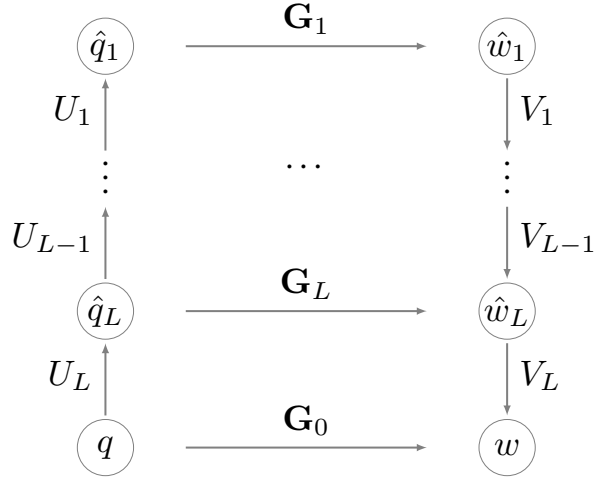

\subsection{Shaping the FMMnet}

Having a fixed sparsity pattern obtained on the first step,
we substitute linear mappings with more general non-linear maps
structured as artificial neural networks. During learning, we assume
that coordinates are known and the cluster construction procedure
is fixed. For a linear transform, we store linear transformations
between the parents and the children of block cluster trees and also linear transformations between the nodes. However, the latter ones are to be replaced by nonlinear mappings in case of FMMNet.
The simplest way to do this is to build feedforward neural networks with block sparse linear parts with fixed sparsity patterns. 
Thus the equation~\ref{eq:h2} takes the following form:
\begin{equation}
    \label{eq:clbl}
    w = {\bf G}_0 \left[ q \right] + \sum_{l=L}^{1} \left( U_L \ldots U_{L-l+1}
    \right) {\bf G}_l \left[ \left( V_{L-l+1} \ldots V_{L} \right) q \right], \\
\end{equation}
where ${\bf G}_l \left[ z \right] $, $l=0\ldots L$ are neural networks with a certain structure which in general may be arbitrary. For example, if one uses
fully-connected layer as an elementary unit of such a network, it can be expressed as

\begin{equation}
{\bf G}_l[z] = Q_{\zeta} \left[ Q_{\zeta-1} \left[ \ldots Q_{1} \left[ z \right]
\right] \right], \quad Q_i \left[t \right] = f_i \left[ W_i t \right],
\end{equation}

where $f_i(t)$ is predefined nonlinearity (for instance, hyperbolic tangent or rectified linear unit), $W_i$ defines linear transformation for $i$-th layer.

Our idea is to utilize a sparsity pattern of original $M_l$ coming from a certain
$\mathcal{H}^2$~matrix. Block sparse structure of these matrices allows us to
split each of them into a set of non-zero blocks supplied by their own
nonlinearities:

\begin{equation}
\left( Q_{i} \left[t \right] \right)_{k} = f_i \left[ \sum\limits_{l \in L_i}
h_{kl} \left[ W_{k l} t_{l} \right] \right],
\end{equation}

where $k$ is a number of output vector block, $l$ - number of input vector block,
$h_{kl} \left[ z \right]$ and $f_i \left[ z \right]$ - nonlinearities,
$W_{kl}$ - blocks of matrix, certain number of which are zero-valued.

Another way to introduce non-linearity into the model is to replace
each non-zero block of block sparse matrices $M_l$ and $C_l$ with its own non-linear processing unit (e.g., feedforward network).

Therefore, as generalization, we constructed a block sparse layer
which preserves the given sparsity template, and it was used as a basic unit for
${\bf G}_l \left[ z \right]$ in our design. We used rectified linear unit (ReLU)
as nonlinearity because it is close to a linear transformation and preserves
properties that make the model easy to optimize and generalize well
\citep{Goodfellow-et-al-2016}, though any other nonlinear function suited for
the learning process may be selected.

It should be noted that we have assumed that the sizes of parameters were defined
a priori. As consequence, we propose to use a certain pre-computed
$\mathcal{H}^2$~matrix structure ($\mathcal{H}^2$-skeleton) as supplier of
related hyperparameters. In practice we propose to use the following strategy:
    \begin{itemize}
        \item If node sets $x$ and $y$ are given for $q(x)$ and $w(y)$,
        then repeat $\mathcal{H}^2$~building procedure
        (\ref{sec:hs,sec:cafb,sec:llc,sec:lc}).
        One should evade to store matrices $C, U_l, V_l, M_l, \; l=1\dots L$
        explicitly, because only their sizes are to be used further. We refer
        to the sequence of such sizes as {\bf $\mathcal{H}^2$~skeleton}.
        \item In other case assume that sets $x$ and $y$ are located on tensor
        grid of the appropriate size and dimension and build
        $\mathcal{H}^2$~skeleton for it.
    \end{itemize}

\subsection{Comparison to other architectures based on $\mathcal{H}^2$-matrix}

Alternative implementation of $\mathcal{H}^2$~based neural network is presented
in~\citet{fan-h2mix-2019}. Authors explored a similar idea to consider
$\mathcal{H}^2$~matrix by vector product as a forward pass of a certain neural
network, but the main difference consists in generalization of $C$ and $M_l$
factors. As basic units authors used either
locally-connected or convolutional layer. The resulting
networks were called as MNN-$\mathcal{H}^2$-LC and MNN-$\mathcal{H}^2$-mix
respectively. However, due to the fact that both locally-connected and
convolutional layers are the representations of banded matrices, the
MNN-$\mathcal{H}^2$-LC and MNN-$\mathcal{H}^2$-Mix networks are
closer to HSS-matrix based architecture rather than the $\mathcal{H}^2$~matrix
one.
\begin{rmk}
    Hierarchically semiseparable (HSS) matrices~\citep{ChLyons-hss-2005,martinsson-hss_integr-2005}
    are the one dimensional versions of $\mathcal{H}^2$~matrices.
\end{rmk}

Thus MNN-$\mathcal{H}^2$-LC network inherits disadvantages of HSS matrices
such as memory issues for problems with a significant off-diagonal part in the
close matrix. Usage of convolutional layers solved the memory problem in case
of MNN-$\mathcal{H}^2$-Mix.

In addition, in Section 4 we provide an example of operator that is difficult to
be approximated by network of such architecture.

\section{Computational experiments}
\label{sec:exp}

\subsection{Implementation details}
All models and computational experiments were implemented using Python
programming language under Anaconda distribution \citep{anaconda}, which
includes various pre-built packages for scientific computing. In this study
the following packages were used:
numpy \citep{oliphant2006guide},
scipy \citep{scipy},
matplotlib \citep{hunter2007matplotlib}.
The $\mathcal{H}^2$-NN model implementation is based on PyTorch package
\citep{pytorch}; for experiments with MNN-H2 network the source code provided
by authors of paper~\citet{fan-hmix-2018} was used which is available at
\url{https://github.com/ywfan/mnn-H2} and based on Keras/TensorFlow frameworks
\citep{chollet2015keras,tensorflow2015-whitepaper}. Some computational experiments
are structured as Jupyter Notebooks \citep{Kluyver:2016aa}.

\subsection{Radiative transfer equation (RTE)}
\label{sec:rte1}

To compare $\mathcal{H}^2$-NN with models proposed in~\citet{fan-h2mix-2019}, we
considered the one-dimensional Radiative transfer equation~\citep{chand-rte-2013}
(RTE):

\begin{equation}
    \label{eq:dif_rte}
    \begin{matrix*}[l]
    v \nabla_x \varphi (x,v) + \mu_t(x,v) = \mu_s(x)u(x)+f(x), \; \mathrm{~in~}
    \Omega \times \mathcal{S}^{d-1}, \Omega \in \mathbb{R}^d, \\
    \varphi(x,v) = 0, \; \mathrm{~on~} \{(x,v)\in \partial\Omega\times
    \mathcal{S}^{d-1}:n(x)\dot v < 0\},\\
    u(x) = \frac{1}{4\pi} \int_{\mathcal{S}^{d-1}}\varphi (x,v)dv
    \end{matrix*}
\end{equation}

Authors of~\citet{fan-hmix-2018} proposed the way to come from differential
form~\eqref{eq:dif_rte} into the integral one,
\begin{equation}
u(x) = \mathcal{I}(\mu_s(x)), \quad \text{or } \mathcal{I}: \, \mu_s(x)
\rightarrow u(x)
\end{equation}

Having a dataset of paired discretized originals and images
$\left(\mu_s(x_h)^{(i)}, u(x_h)^{(i)}\right)$, $i=1\ldots N_{\text{s}}$, one
may learn the integral mapping using a certain predefined model, for example,
parametric model $g(z, \theta)$ with trainable parameters $\theta$. Training
process is guided by optimizing specified loss functional that measures
quantitatively how output of model with current values of parameters differs
from the ground truth:

\begin{equation}
\widehat{\theta} = \arg \min\limits_{\theta}\mathcal{L}\left( u(x_h),
g\left( \mu_s(x_h), \theta \right) \right) =
\arg \min\limits_{\theta} \frac{1}{N} \sum\limits_{i=1}^{N}
\mathcal{L}_i\left( u(x_h)^{(i)}, g\left( \mu_s(x_h)^{(i)}, \theta \right)
\right)
\end{equation}

Because both inputs and outputs in this task are real-valued vectors, one may
use any loss functional valid for regression-like problems. In this work the
squared relative residual was selected in this regard:
\begin{equation}
\mathcal{L}_i\left( u(x_h)^{(i)}, g\left( \mu_s(x_h)^{(i)}, \theta \right)
\right) = \frac{\| u(x_h)^{(i)} - g\left( \mu_s(x_h)^{(i)}, \theta \right)
\|_F^2}{\| u(x_h)^{(i)} \|_F^2}
\end{equation}

Dataset for this problem was generated using the source code provided by the
authors of paper~\citet{fan-hmix-2018} (see Section 4.1).
All samples were computed on the coinciding sets of sources $x$ and receivers
$y$, $x = y$. For three variants of grid sizes, $N = 320, 640, 1280$
equal number of samples $N_{\text{s}} = 20000$ were generated. To control
generalization ability of each model, holdout validation scheme was used with
splitting dataset into training and validation parts in the ratio of $2:1$.

In addition to model proposed in this study, three other ones were selected
for comparison, namely simple convolutional model (conv),
MNN-$\mathcal{H}^2$-Mix and MNN-$\mathcal{H}^2$-LC
(both from~\citet{fan-hmix-2018}). All models were trained using the following
hyperparameters:

\begin{itemize}
    \label{it:model_cond}
    \item Optimizer: Adam \citep{kingma-adam-2014}
    \begin{itemize}
        \item Learning rate: lr = 0.0025,
        \item Coefficients used for computing running averages of gradient and
        its square: $\beta=(0.9, 0.999)$,
        \item Term added to the denominator to improve numerical stability
        $\epsilon=10^{-5}$,
    \end{itemize}
    \item Number of layers: $\zeta = 3$,
    \item Nonlinearity: ReLU,
    \item Number of iterations: $N_{\mathrm{it}} = 2000$,
    \item Model initialization: Glorot uniform initializer,
          also known as Xavier uniform initializer~\citep{glorot-init-2010}.
          
\end{itemize}

\begin{figure}[H]
\begin{subfigure}{.33\linewidth}
\centering
\includegraphics[width=1.\linewidth]{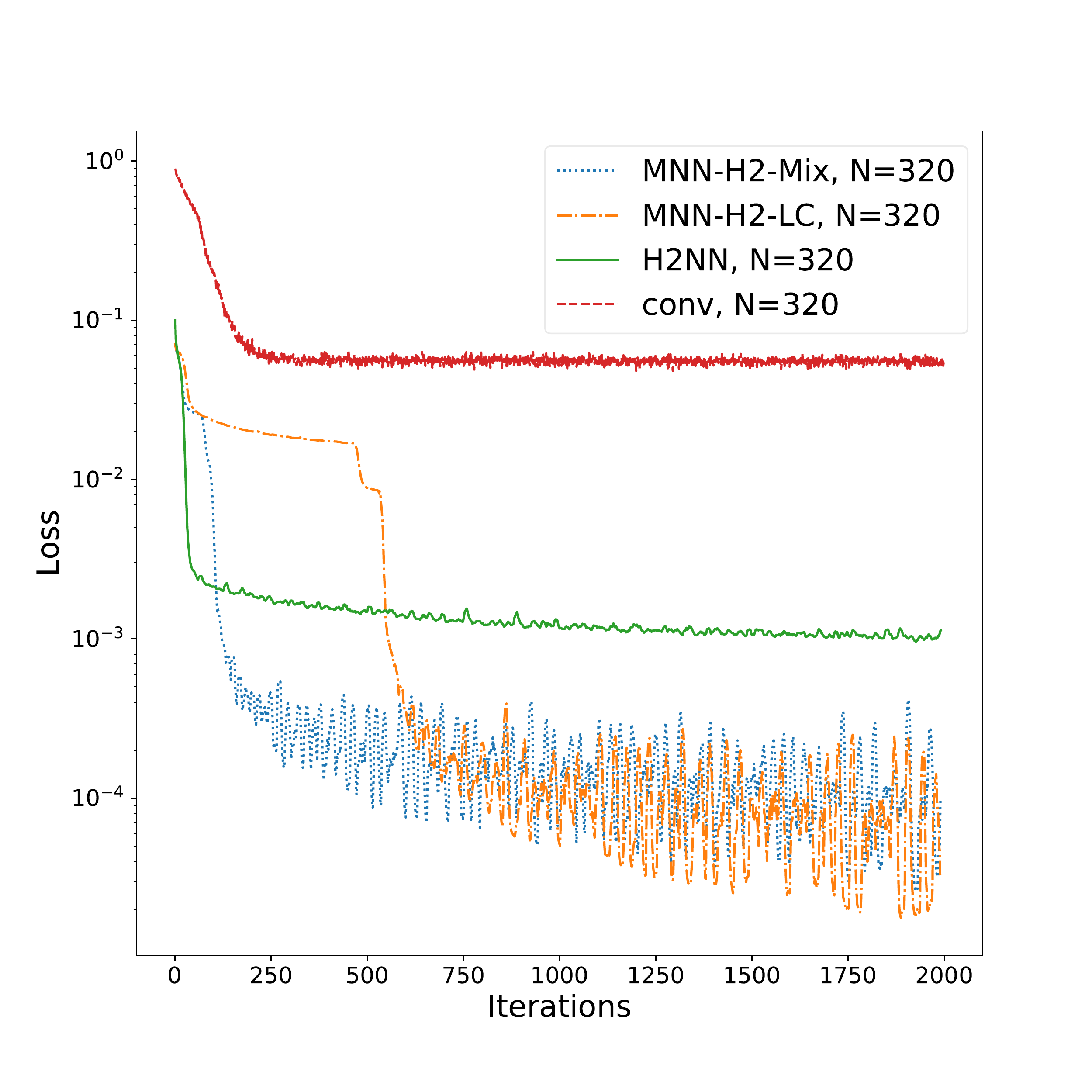}
\caption{$\,$}
\label{fig:sub1}
\end{subfigure}%
\begin{subfigure}{.33\linewidth}
\centering
\includegraphics[width=1.\linewidth]{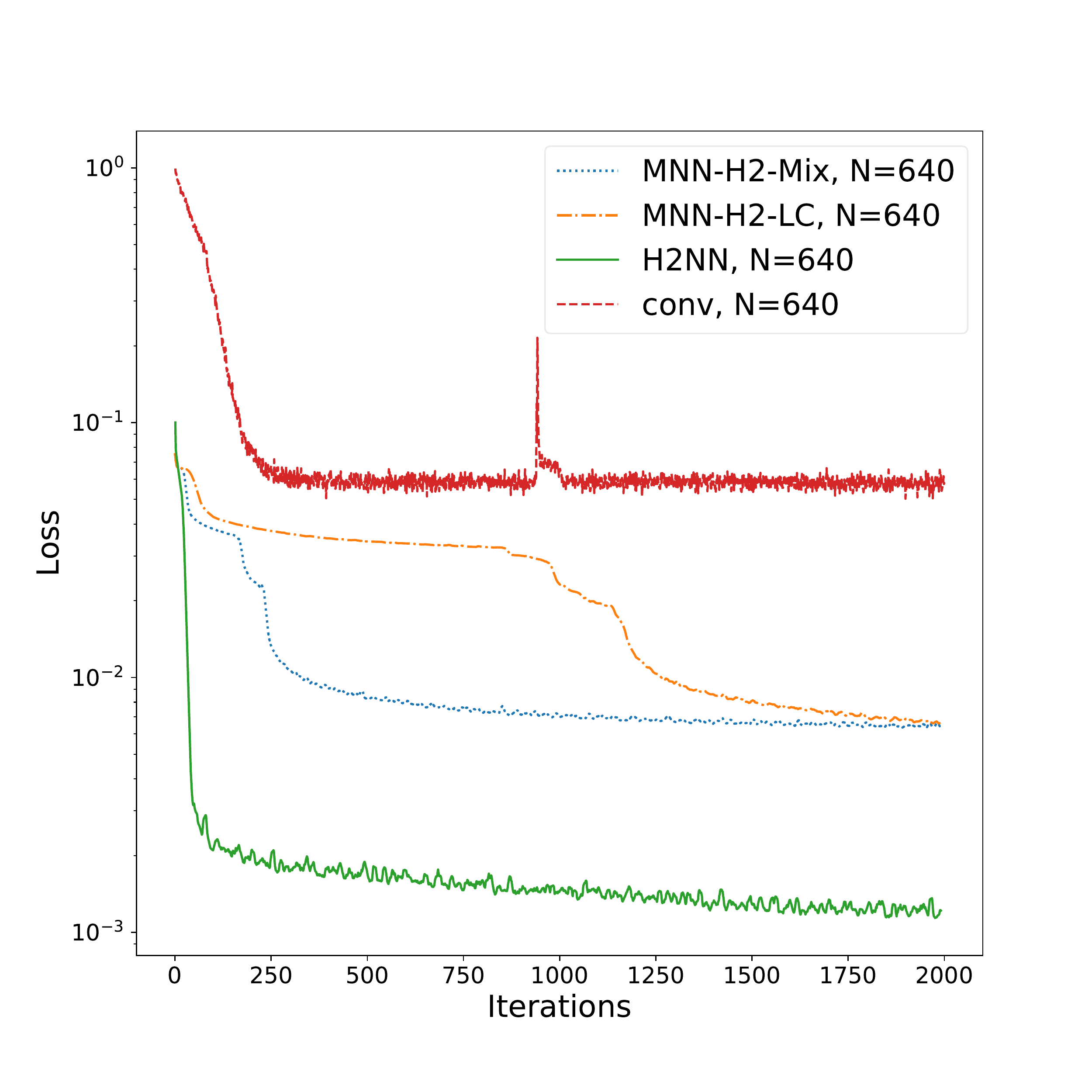}
\caption{$\,$}
\label{fig:sub2}
\end{subfigure}
\begin{subfigure}{.33\linewidth}
\centering
\includegraphics[width=1.\linewidth]{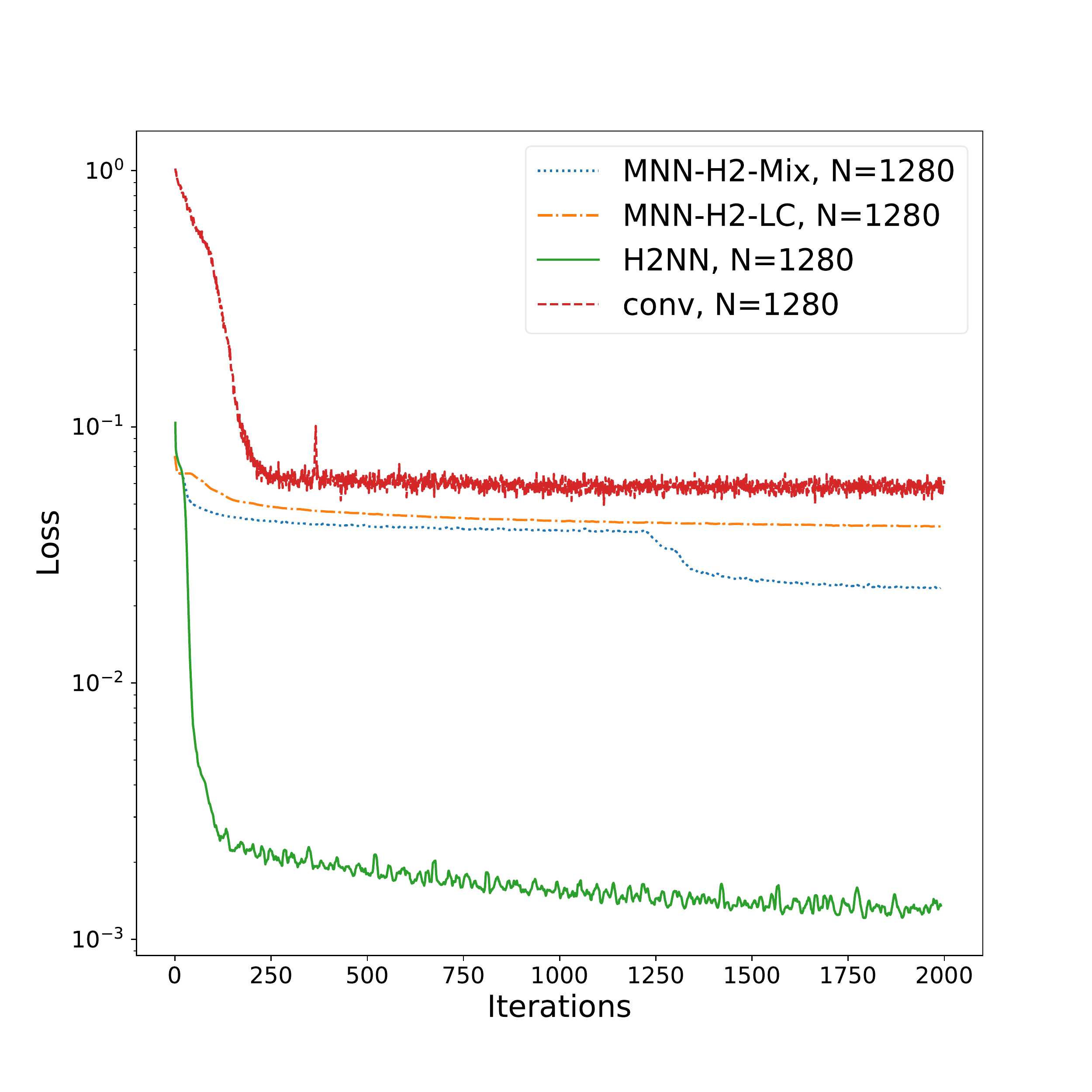}
\caption{$\,$}
\label{fig:sub3}
\end{subfigure}
\caption{Loss decrease during optimization process for RTE 1D with different
grid sizes $N$: (a) $N=320$; (b) $N=640$; (c) $N=1280$. Curves legend:
dotted blue line for multiscale neural network based on hierarchical nested
bases with convolutional and locally-connected layers (MNN-H2-Mix); dash-dotted
orange line for multiscale neural network based on hierarchical nested bases
with locally-connected layers (MNN-H2-LC); solid green line for
$\mathcal{H}^2$-based neural network (H2NN, proposed approach); dashed red for
simple convolutional network (conv).}
\label{fig:rte_loss}
\end{figure}

Figure~\ref{fig:rte_loss} displays learning process it terms of evoluting
loss functional. Average values of time per iteration presented in
Table~\ref{tab:rte_it_time}.

In comparison to proposed network, both MNN-$\mathcal{H}^2$-Mix and
MNN-$\mathcal{H}^2$-LC converged faster to the lower values of loss functional
in case of $N=320$. But for higher dimensionality these models saturated to
larger loss values while proposed model quickly decreased
to lower ones. As expected, simple but not specific network showed the worst
performance for all considered dimensions. These observations suggest
evidence that $\mathcal{H}^2$-NN model is more scalable.
Table~\ref{tab:rte_trte} contains information about train and validation
errors (mean relative residuals) after $N_{\mathrm{it}} = 2000$ iterations.

\begin{table}[H]
\centering
\scalebox{1.}{
    \input{figs/rte_table.tex}
}
\caption{RTE 1D, comparison of mean relative residual errors for different
methods and grid sizes $N=320,640,1280$. Each cell is formatted as
$\varepsilon_{\text{train}}/\varepsilon_{\text{test}}$, where
$\varepsilon_{\text{train}}$ and $\varepsilon_{\text{test}}$ are mean relative
residual errors measured on training and test sets.
}
\label{tab:rte_trte}
\end{table}

\begin{figure}[H]
\centering
\begin{subfigure}[t]{0.49\textwidth}
    \centering
    \includegraphics[width=1.\linewidth]{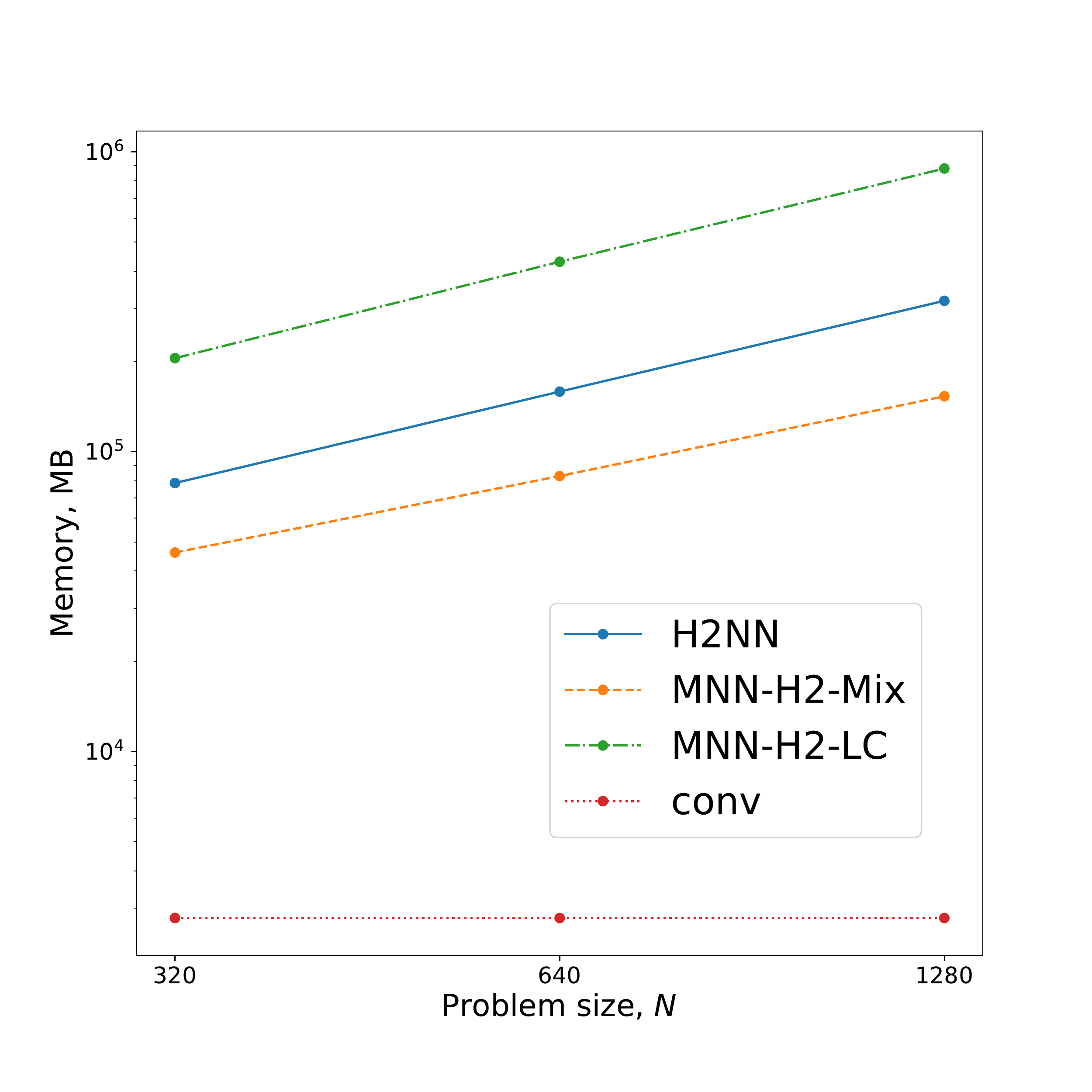}
    \caption{$\,$}
\end{subfigure}
\begin{subfigure}[t]{0.49\textwidth}
    \centering
    \includegraphics[width=1.\linewidth]{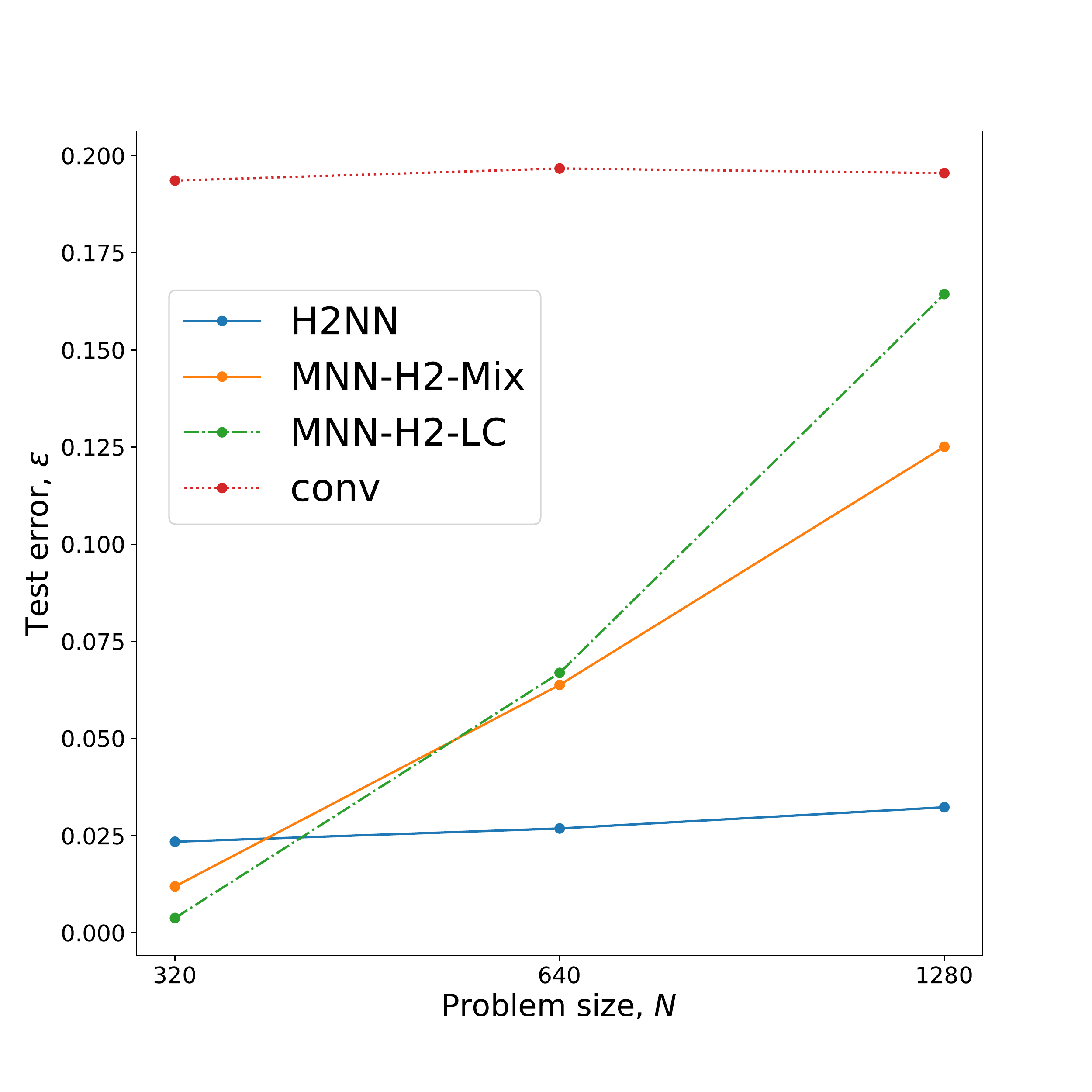}
    \caption{$\,$}
\end{subfigure}
\caption{Memory and error comparison of $\mathcal{H}^2$-NN,
MNN-$\mathcal{H}^2$-Mix and MNN-$\mathcal{H}^2$-LC models. (a) memory required for the model storage; (b)
test error}
\label{fig:rte1d_mem_err}
\end{figure}

Mean relative residuals measured on validation part of datasets as well as memory
consumption of considered models are shown in Figure~\ref{fig:rte1d_mem_err}.

According to Figure~\ref{fig:rte1d_mem_err}, proposed architecture
($\mathcal{H}^2$-NN) has the best scalability with moderate memory consumption,
surpassed only by MNN-$\mathcal{H}^2$-Mix model.

\subsection{Custom operator example}
\label{sec:ce}

It is known that convolution operation may be represented as multiplication by
Toeplitz matrix. Experience shows \citep{fan-h2mix-2019} that
$\mathcal{H}^2$-inspired architectures with convolutional layers
(MNN-$\mathcal{H}^2$-Mix) are good approximators for matrices structured closely
to Toeplitz ones. A natural question is how well it capture more general
dependencies. We considered the following weights:
\begin{equation}
a_{ij} = \frac{1}{\| x_i - y_j\|} f(i,j),
\end{equation}
where $f(i,j) = \frac{(x_i+y_j)^{2}}{N^{5}}$.

\begin{figure}[H]
\begin{subfigure}{.33\linewidth}
\centering
\includegraphics[width=1.\linewidth]{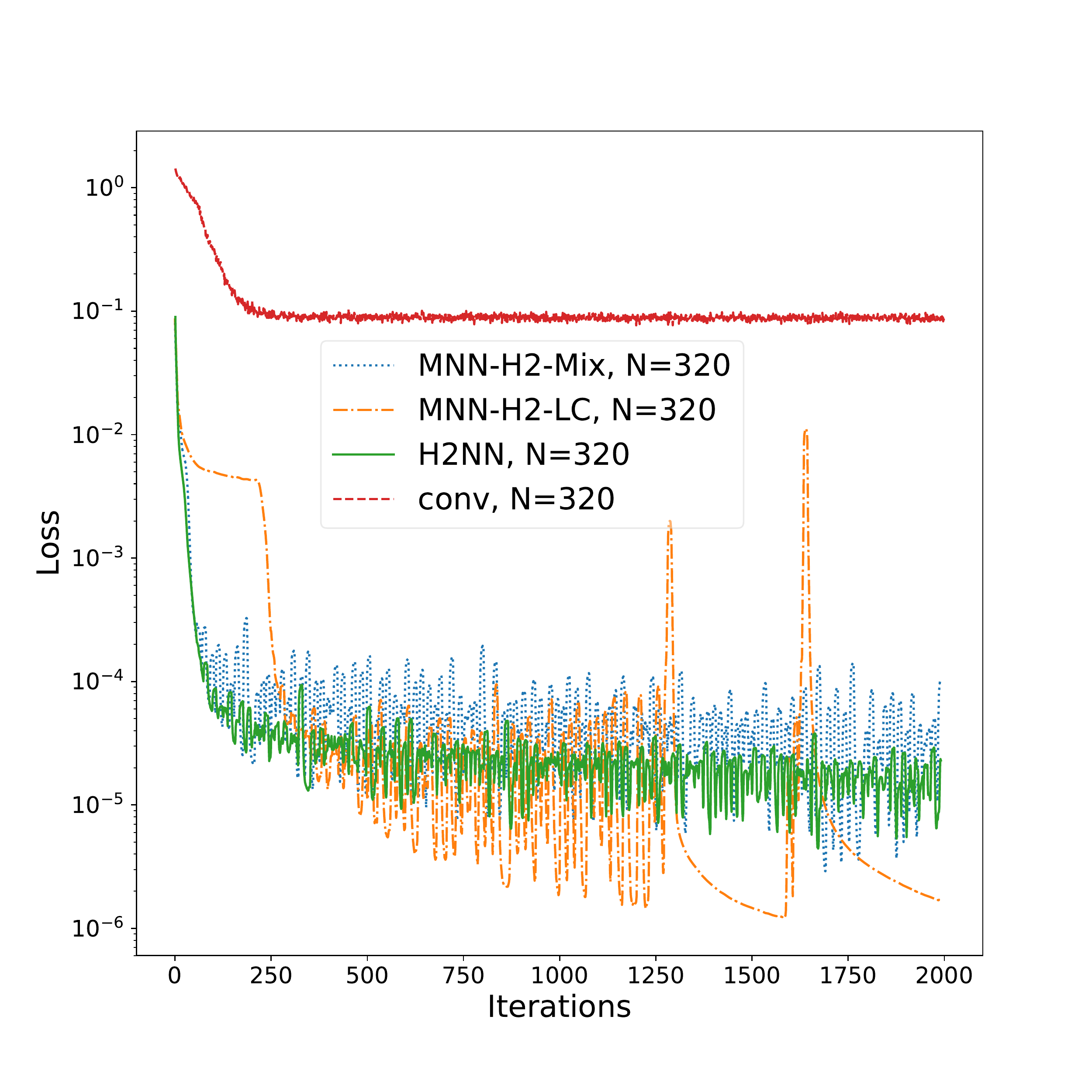}
\caption{$\,$}
\label{fig:sub1}
\end{subfigure}%
\begin{subfigure}{.33\linewidth}
\centering
\includegraphics[width=1.\linewidth]{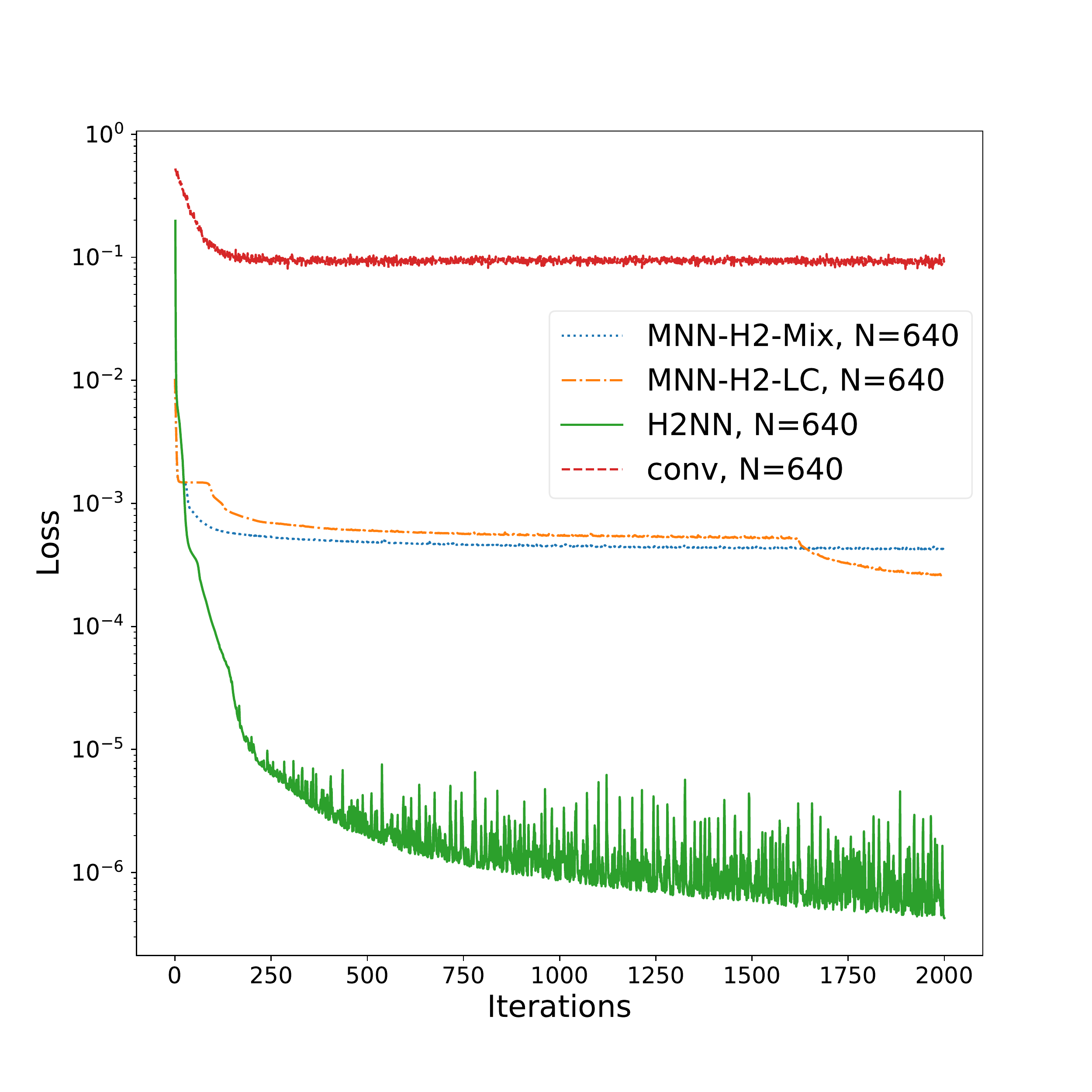}
\caption{$\,$}
\label{fig:sub2}
\end{subfigure}
\begin{subfigure}{.33\linewidth}
\centering
\includegraphics[width=1.\linewidth]{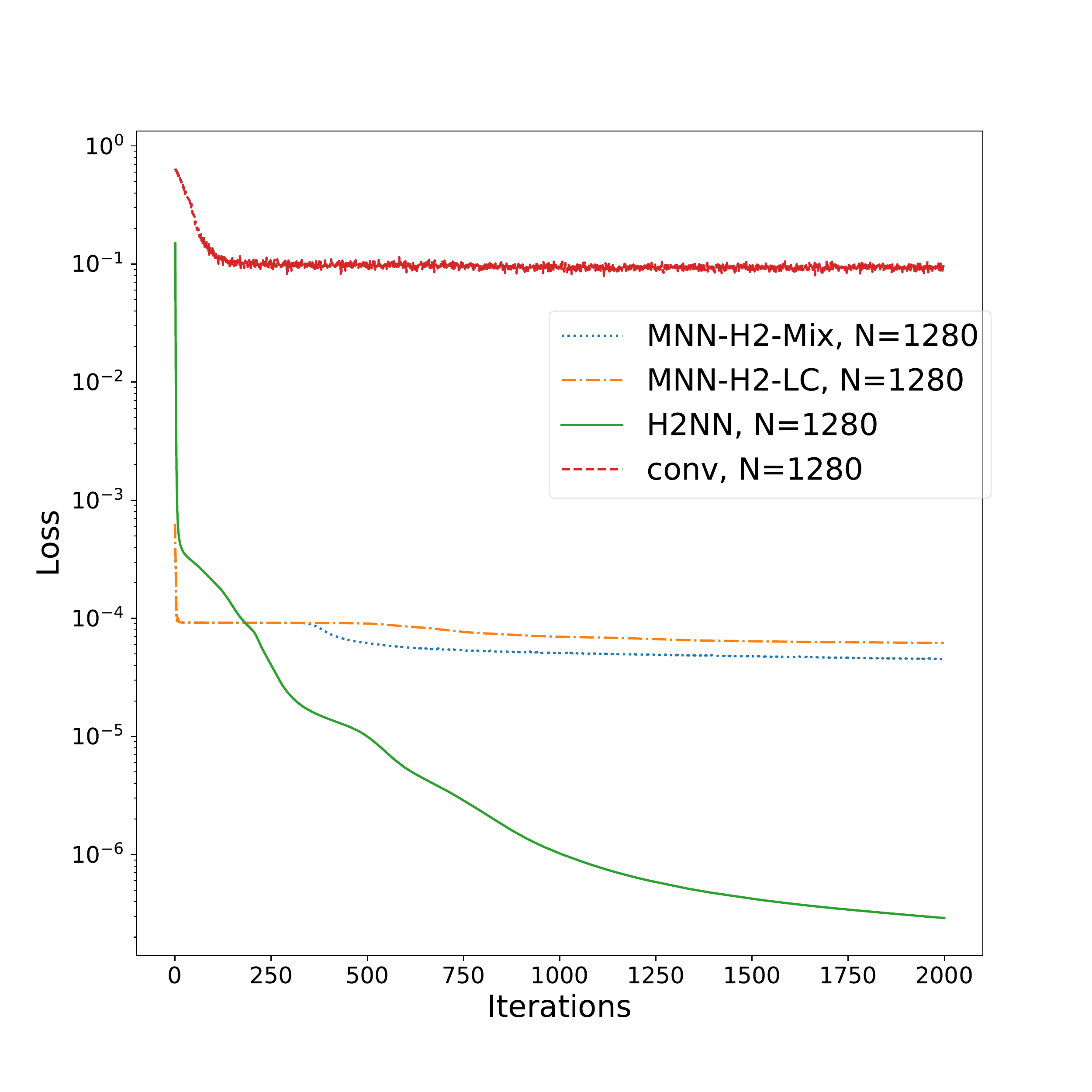}
\caption{$\,$}
\label{fig:sub3}
\end{subfigure}
\caption{Loss decrease during optimization process for custom operator example
 with different grid sizes $N$: (a) $N=320$; (b) $N=640$; (c) $N=1280$. Curves
 legend: dotted blue line for multiscale neural network based on hierarchical
 nested bases with convolutional and locally-connected layers (MNN-H2-Mix);
 dash-dotted orange line for multiscale neural network based on hierarchical
 nested bases with locally-connected layers (MNN-H2-LC); solid green line for
 $\mathcal{H}^2$-based neural network (H2-NN, proposed approach); dashed
 red for simple convolutional network (conv).}
\label{fig:ce_loss}
\end{figure}

\begin{figure}[H]
    \centering
\begin{subfigure}[t]{0.49\textwidth}
    \centering
        \includegraphics[width=1.\linewidth]{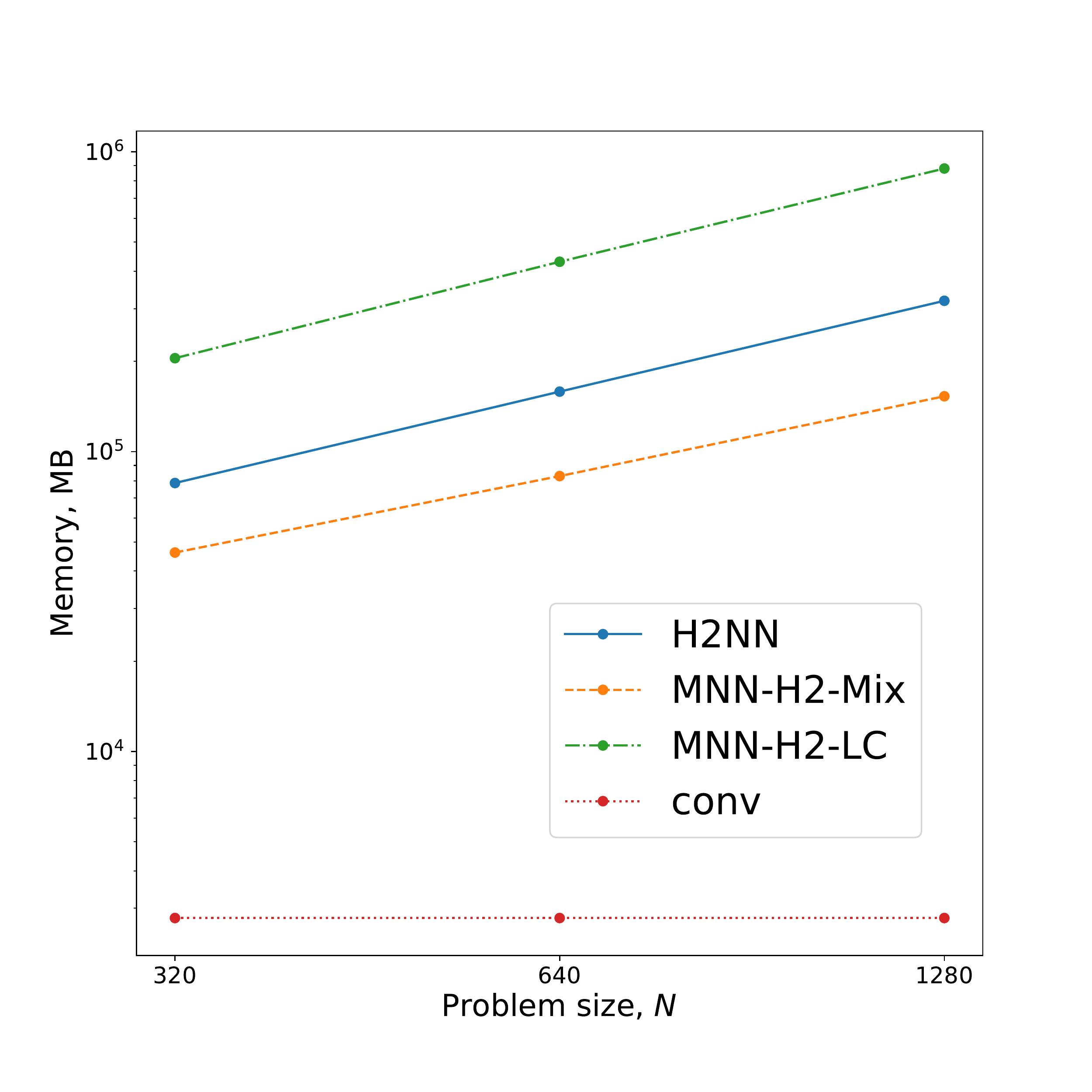}
        \caption{$\,$}
\end{subfigure}
\begin{subfigure}[t]{0.49\textwidth}
    \centering
    \includegraphics[width=1.\linewidth]{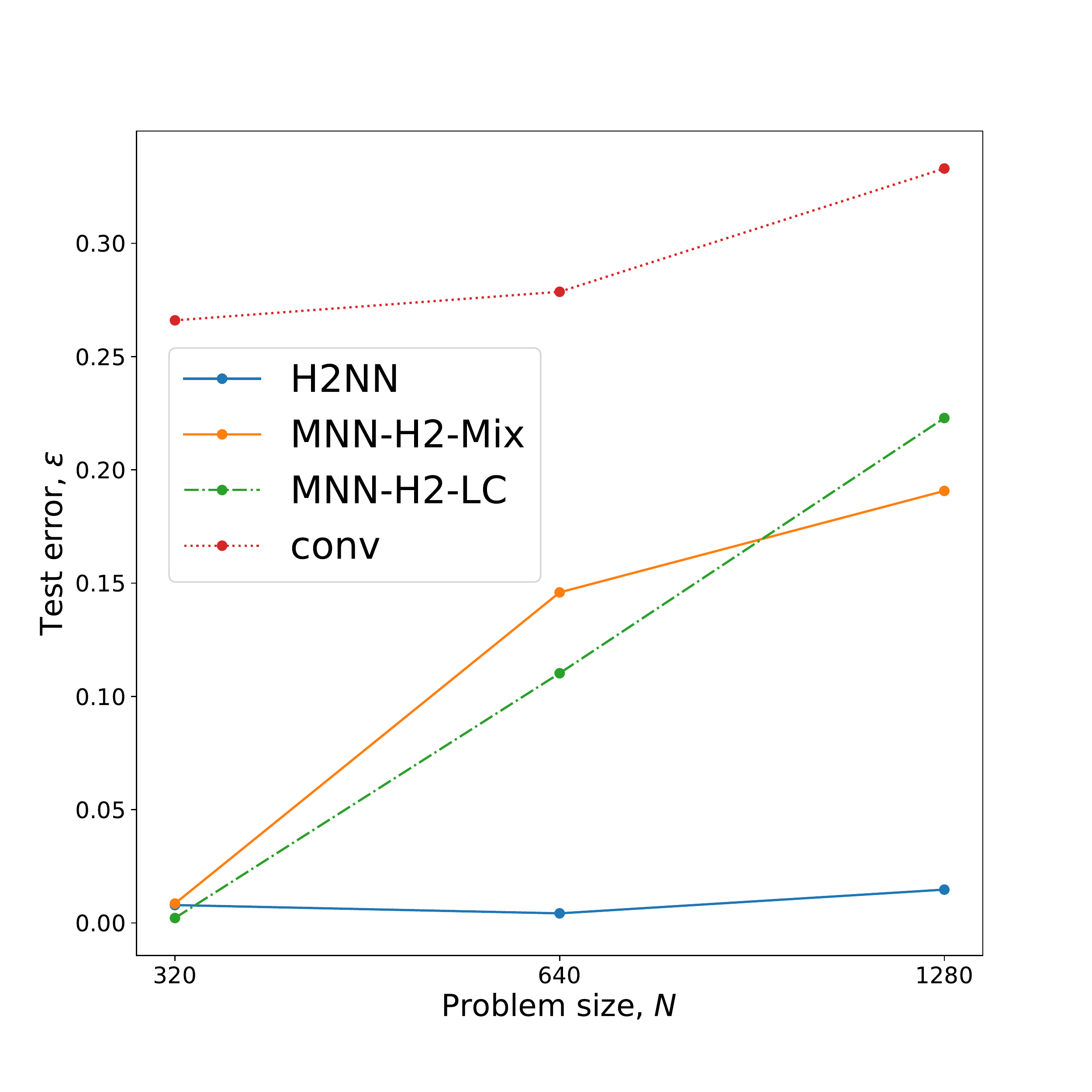}
    \caption{$\,$}
\end{subfigure}
\caption{Memory and error comparison of $\mathcal{H}^2$-NN,
MNN-$\mathcal{H}^2$-Mix and MNN-$\mathcal{H}^2$-LC models. (a) memory required for the model storage; (b)
test error}
\label{fig:ce_mem_err}
\end{figure}

Inputs were similar to ones used in section~\ref{sec:rte1}, and outputs were
generated by this custom weighting procedure. The same setup as in previous
section was used to train and validate conv, MNN-$\mathcal{H}^2$-Mix/LC and
$\mathcal{H}^2$-net models. Training procedure is visualised in
Figure~\ref{fig:ce_loss}. Average values of time per iteration presented in
Table~\ref{tab0001}. Mean residual errors and memory consumption of
models are shown in Figure~\ref{fig:ce_mem_err}.

\begin{table}[H]
\centering
\scalebox{1.}{
    \input{figs/ce_table.tex}
}
\label{tab0001}
\caption{Custom operator example, comparison of mean relative residual errors
for different methods and grid sizes $N=320,640,1280$. Each cell is
formatted as $\varepsilon_{\text{train}}/\varepsilon_{\text{test}}$,
where $\varepsilon_{\text{train}}$ and $\varepsilon_{\text{test}}$ are mean
relative residual errors measured on training and test sets.}
\label{tab:ce_trte}
\end{table}

MNN-$\mathcal{H}^2$-Mix showed nearly comparable with MNN-$\mathcal{H}^2$-LC
model behaviour, though had better performance for lower problem
size, $N=320$ (Figure~\ref{fig:ce_loss}, (a); Table~\ref{tab:ce_trte}). However,
both models were outperformed by $\mathcal{H}^2$-NN model in higher dimensions,
$N=640, 1280$ (Figure~\ref{fig:ce_loss}, (b), (c); Table~\ref{tab:ce_trte}).
Memory consumption preserved the trend outlined in the previous section.

\section{Related work}
\label{sec-rel}

Structured low-rank matrix formats are actively used to solve partial
differential equations (PDEs) and integral equations (IEs)
~\citep{GrRo-fmm-1987,tee-mosaic-1996,bebendorf-existence-2003,white-schur-2009}.
Some of them were intentionally designed to reduce memory and computational
costs. Being the particular example of such formats, $\mathcal{H}$-matrices
~\citep{hackbusch-h-1999} together with mosaic-skeleton ones
~\citep{tee-mosaic-1996} allow to store dense matrices by $\mathcal{O}(N \log N)$
memory cost and to perform matrix-vector multiplication by
$\mathcal{O}(N \log N)$ operations. $\mathcal{H}^2$~matrices
~\citep{hackbusch-h2-2000,Borm-h2-2010,Bebendorf-h2-2012}
which are the central point of the fast multipole method (FMM) and additionally reduce
these costs to $\mathcal{O}(N)$. At the same time, solving the linear
system with $\mathcal{H}$ and $\mathcal{H}^2$~matrix is challenging.
The complicated structure of block low-rank matrices makes it hard to build a
direct solver. Another issue regarding the direct solution is preserving the
$\mathcal{O}(N)$ memory complexity during computations. Despite the theoretical
knowledge that $\mathcal{H}^2$~matrices (with certain conditions) have
$\mathcal{H}^2$~matrix as an inverse
~\citep{bebendorf-existence-2003,Beb-hlu-2005}, the implementation of the
$\mathcal{H}^2$~direct solver is a work in progress.

The rapidly growing area of artificial neural networks allowed to propose differently
approaches to the solution of PDEs
\citep{lagaris-nn_pde-1998,chiaramonte-nn_pde-2013,baymani-nn_pde-2010} and
IEs \citep{elshafiey-nn_ie-1991,vemuri-nn_ie-1991,effati-nn_ie-2012}. In general,
the main idea is to teach a neural network to act as a transformation of the
right-hand side to the solution. Theoretical results state that neural networks
are universal approximators of functional dependencies: according to Cybenko
theorem (universal approximation theorem)~\citep{hornik-uat-1989,
hornik-aut-1991}, feed-forward network with 2 layers (one hidden layer)
of finite sizes and mildly restricted nonlinearity is a universal approximator
for continuous functions on compact sets in $\mathbb{R}^d$. However, 
in many practical applications such as speech recognition, image segmentation,
generative modelling, multi-layered as well as non-feedforward architectures
are actively used.

Typical architecture used for solving PDE and IE is a feedforward network
with a single hidden layer~\citep{chiaramonte-nn_pde-2013,baymani-nn_pde-2010,
elshafiey-nn_ie-1991,vemuri-nn_ie-1991,piscopo-nn_pde-2019,effati-nn_ie-2012}
or with multi-layered structure~\citep{ramchoun-mlp_pde-2016,
lagaris-nn_pde-1998}. These models have fully-connected layers with dense
unstructured matrices leading to memory issues.

The idea of construction $\mathcal{H}$-, $\mathcal{H}^2$-based neural networks
is not unique, and several implementations of it are known
\citet{fan-hmix-2018,fan-h2mix-2019} where authors used band matrix for close
and tree-to-tree transfer matrices. This paper develops the idea of
$\mathcal{H}^2$-NN architecture, proposing a new method of layers construction
and showing on practice benefits of this method.

\section{Conclusions}
\label{sec-con}

In this work, we proposed a new $\mathcal{H}^2$-based architecture and showed its
benefits over the existing analogs.

In the context of partial differential and integral equations used for
modeling physical processes, machine learning holds the promise of being
able to capture relationships between observable measurements in complicated
conditions (e.g., complex geometry). Moreover, trainable parametric models may
potentially generalize to various input conditions. In that case, neural network
architectures inspired by established constructions like
$\mathcal{H}^2$-matrices would be easier to analyze.

It is worth noting that a common problem for every neural network-based
method for solving PDEs and IEs is a large learning time. However, as soon as
the network is trained, the solution is to be computed quickly by one forward
pass. Since training is required at once, such an approach is
prospective in case of a large number of typical problems with the shared operator
and different right-hand sides.

In future work, we will continue researching $\mathcal{H}^2$-networks generated by
block-sparsity patterns emerging in $\mathcal{H}^2$-matrices. Complementary to
problems that arise in mathematical physics, $\mathcal{H}^2$-networks seem to be
applicable in conventional machine learning problems.

\bibliography{iclr2020_conference}
\bibliographystyle{iclr2020_conference}

\appendix

\section{Appendix: Basic concepts of $\mathcal{H}^2$~matrix}
In general, computation of matrix by vector product requires
$\mathcal{O}(N^2)$ operations. and storage of corresponding matrix requires $\mathcal{O}(N^2)$ items of memory. However, additional knowledge about the problem may reduce the computation and storage complexity to $\mathcal{O}(N)$.

\subsection{Separation property}
\label{sec:sp}

Low-rank matrices is an example of structured matrices that have reduced number
of parameters, $\mathcal{O}(2NR)$, where $R$ is a matrix rank, and $R < N$. Due
to its compressing property, low-rank assumption is widely used in many
applications of computational science. 
Matrix $A$ from equation~\ref{eq:matrix-wsum} also shares the low-rank constraint
if sets $x$
and $y$ are {\bf spatially separated}. We refer to this assumption as
separation property of the problem.

\begin{rmk}
    \label{rmk:sp}
Spatial separation implies that the distance $\rho(x_c, y_c)$ between centers
$x_c$ and $y_c$ of the bounding boxes of sets $x$ and $y$ is larger than some
constant $\varrho$ scaled by the size of bounding boxes $\alpha_{\text{bb}}$.
\begin{equation}
\rho(x_c, \, y_c) > \varrho\alpha_{\text{bb}}
\end{equation}
\end{rmk}

The assumptions is valid for N-body problem as well as for IEs with smooth
kernels \cite{greengard-fmm-1987,greengard-fmm-1998}.

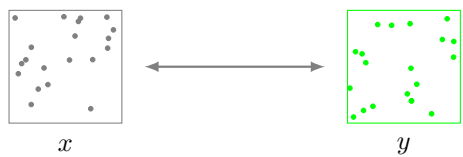
\begin{figure}[H]
\centering
\input{figs/separation_property.tex}
\caption{Example of sets that share separation property: two sets $x$ and $y$
that have bounding boxes distant enough from each other}
\label{fig:sp}
 \end{figure}

\subsection{Hierarchical grid}
\label{sec:hs}

In general case sets $x$ and $y$ are not spatially separated in a straight way.
Moreover, in many problems these sets are coincident, thus a separation
property does not hold for matrix $A$ directly. However, one may consider
subsets of $x$ and $y$, and for some pairs of these subsets the spatial
separation property is valid. Consequently, there are submatrices (blocks) of
matrix $A$ with low-rank. That is the basic idea behind block low-rank matrices.

The natural way to obtain subsets is to split sets $x$ and $y$ into $\eta$
sub-blocks {\bf hierarchically}. Sets $x$ and $y$ are separated into  $\eta$
arbitrary equal subsets~\citep{williams-inertial_bisection-1994}. Then each
subset is recursively separated into its $\eta$ subsets until the stopping
criteria is reached.

\begin{figure}[H]
\centering
\begin{subfigure}[t]{0.45\textwidth}
    \centering
    \input{figs/grid_a.tex}
    \caption{$\,$}
    \label{fig:hier}
\end{subfigure}
\begin{subfigure}[t]{.45\textwidth}
    \centering
    \input{figs/grid_b.tex}
    \caption{$\,$}
    \label{fig:close}
\end{subfigure}
 \caption{Grid operations in two-dimensional case: (a) hierarchical nested division
 into blocks; (b) close-far assignments to the blocks at the same level.
 Close blocks are highlighted by blue color. The latter operation may be
 performed in various ways.}
 \label{fig:grid}
 \end{figure}
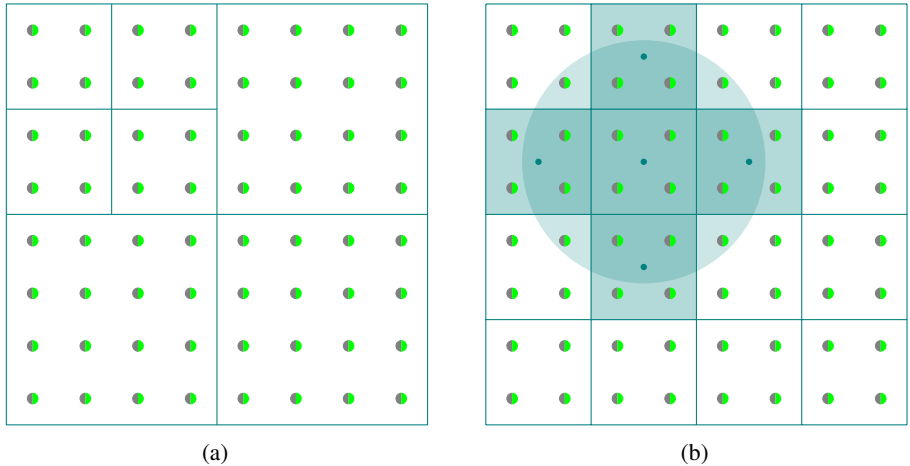

\begin{example}
\label{ex:grid}
As an example consider $x$ and $y$ on the same uniform tensor grid in 2d.
We separate squares into four equal sub-squares ($\eta = 4$), see Figure~\ref{fig:hier}.
The number of boxes on level $l$ of the grid is~$H_l = \eta^{l}$.
\end{example}

\begin{rmk}
    Enumeration of 2d/3d clouds of points is crucial to constitute sets $x$ and
    $y$. Essentially, one should provide a rule for conversion of coordinate
    arrays into flat indices. This operation is equal to the permutation of
    vectors $q$ and $w$, and rows and columns of matrix $A$.
    The standard convention is a depth-first box-wise numeration.
\end{rmk}

\subsection{Close and far blocks}
\label{sec:cafb}

In $\mathcal{H}^2$~matrices, nested grid can be subdivided into regions of two
types at each level. Consider the level with the smallest box size, $l=L$, and
let us look closer at the structure of source and receiver vectors. The former
one contains the following blocks:

$$\tilde{q}_i = q(\tilde{x}_i), \; i \in 1\dots H_L, \\
  \tilde{w}_j = w(\tilde{y}_j), \; j \in 1\dots H_L,$$
where $\tilde{q}_i,\tilde{w}_j\in\mathbb{R}^{B}$, $B$ is block size, $H_L$ is
number of blocks on level $L$. Thus,
$$q = \begin{bmatrix}
    \tilde{q}_1\\
    \vdots\\
    \tilde{q}_{H_L}
\end{bmatrix}, \;
w = \begin{bmatrix}
    \tilde{w}_1\\
    \vdots\\
    \tilde{w}_{H_L}
\end{bmatrix}, \;
A = \begin{bmatrix}
    A_{11}& \dots & A_{1H_L}\\
    \vdots&\ddots&\vdots\\
    A_{H_L1}& \dots & A_{H_L}\\
\end{bmatrix},$$
where $A_{ij}\in\mathbb{R}^{B\times B}$.

For each receiver box $\tilde{y}_i$, $i \in 1\dots H_L$ let us consider set of
sources boxes $\tilde{x}_j$ within the circle of radius $\rho = \varrho\alpha_{L}$,
where $\varrho$ is a predefined constant, $\alpha_{L}$ is a box size on level $L$.

\begin{figure}[H]
\centering
\resizebox{.8\textwidth}{!}{
    \input{figs/mat_A.tex}
}
\caption{Example of $\mathcal{H}^2$~matrix $A = C+F_L$, $L=4$. Far and close
blocks are highlighted by light and dark blue colors respectively. }
\label{fig:mata41}
\end{figure}
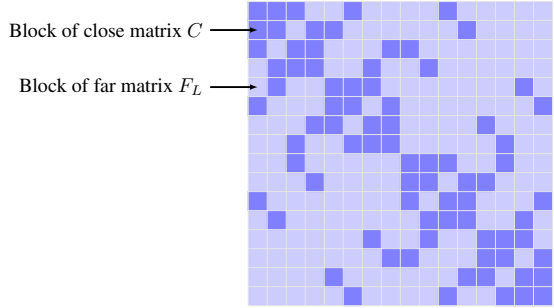

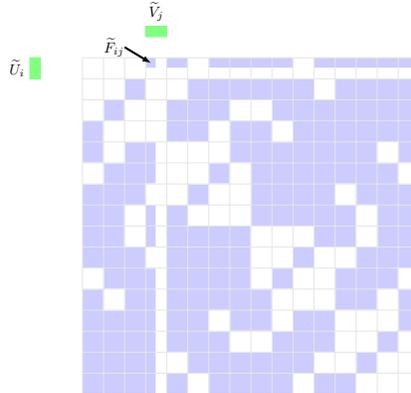
\begin{figure}[H]
\centering
\resizebox{.4\textwidth}{!}{
    \input{figs/bottom_level_compression.tex}
}
\caption{Matrix $F_L$, $L=4$ that contains far blocks of $\mathcal{H}^2$~matrix.
Decomposition of pointed block $F_{ij}$ into $\widetilde{U}_i \widetilde{F}_{ij}
\widetilde{V}_j^T$ product allows to compress it; matrices $\widetilde{U}_i$,
$\widetilde{V}_j$ are assumed to be common for every block of $i$-th row and
$j$-th column respectively.}
\label{fig:mata42}
\end{figure}

In case of receivers close to sources, corresponding blocks $A_{ij}\in
\mathbb{R}^{B\times B}$
of matrix $A$ have {\bf no low-rank} since spatial separation property does not
hold for them. Because they correspond to closer interactions, we call them
{\bf close blocks}.
On the level $l=L$ these close blocks constitute close matrix $C$, see Figure~\ref{fig:mata41}.
Remaining blocks called {\bf far blocks} can be approximated well with low-rank.
The far blocks constitute {\bf far matrix}~$F_L$, see Figure~\ref{fig:mata42}.
Thus the level $L$ structure of $\mathcal{H}^2$~matrix-by-vector product can be
rewritten as
\begin{equation}
w = Aq = Cq + F_Lq
\end{equation}.

 \subsection{Bottom-level compression}
 \label{sec:llc}

Owing to separation property, each far block $F_{ij}\in \mathbb{R}^{B\times B}$
of far matrix $F_L$ has a low rank:
\begin{equation}
F_{ij} \approx \widetilde{U}_i \widetilde{F}_{ij}  \widetilde{V}_j, \quad
\forall i,j \in 1,\dots, H_L,
\end{equation}
where $\widetilde{F}_{ij} \in {\color{black}\mathbb{R}^{r\times r}}$ is a
compressed far block, matrices
$\widetilde{U}_i \in {\color{black}\mathbb{R}^{B\times r}}$ and
$\widetilde{V}_j \in {\color{black}\mathbb{R}^{r\times B}}$ are rectangular
compression matrices.

The key assumption of $\mathcal{H}^2$~matrix is that all the blocks in $i$-th
row have the same left rectangular compression factor~$\widetilde{U}_i$
and all blocks in $j$-th column have the same rectangular compression right
factor $\widetilde{V}_j^{\top}$.
The goal of the compression procedure is to sparsify the matrix $A$ by obtaining
compressed blocks $\widetilde{F}_{ij}$ instead of original blocks $F_{ij}$.

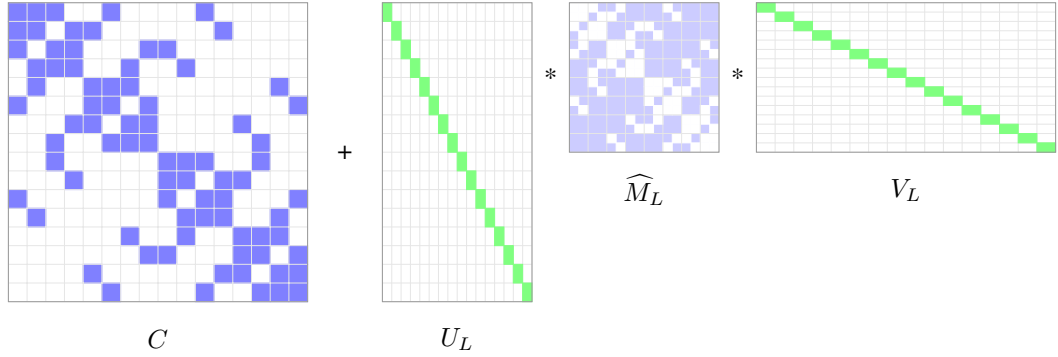
\begin{figure}[H]
\centering
  \input{figs/L_lvl_compr.tex}
  \caption{$\mathcal{H}^2$~matrix, compression at last level, $L$. $C$ matrix
  captures close interaction blocks, $\widehat{M}_L$ is a compressed far
  interaction matrix, $U_L$ and $V_L$ are block diagonal compressing
  matrices.}
  \label{fig:lvl_L_fact}
  \end{figure}

~\\We introduce the block-diagonal rectangular
 compression matrix
 \begin{equation} U_L =
\begin{bmatrix}
\widetilde{U}_1 &0 & 0 \\
0& \ddots &0  \\
0&0 & \widetilde{U}_M
\end{bmatrix},
\label{eq:q}
\end{equation}
where $U_L \in \mathbb{R}^{B H_L\times r H_L}$.
Similarly, for block columns, we obtain the block-diagonal rectangular
 compression matrix \begin{equation} V_L =
\begin{bmatrix}
\widetilde{V}_1 &0 & 0 \\
0& \ddots &0  \\
0&0 & \widetilde{V}_M
\end{bmatrix},
\label{eq:v}
\end{equation}
where $V_L \in \mathbb{R}^{rH_L\times BH_L}$.
For the far matrix, we obtain the factorization
$$F_L = U_L\widehat{M}_LV_L,$$
where $\widehat{M}_L \in \mathbb{R}^{r H_L\times r H_L}$ is a compressed far
matrix. Thus, for the $A$ we obtain
\begin{equation}
    \label{eq:comp_lvl_L}
    w = Cq + U_L\widehat{M}_LV_Lq,
\end{equation}
see Figure~\ref{fig:lvl_L_fact}. 

Consider the matrix $\widehat{M}_L$; our next goal is to compress this matrix.
In order to do that, we move to level $l = L-1$. On this level blocks are
united in groups of $\eta$.

\begin{figure}[H]
\centering
    \resizebox{0.8\textwidth}{!}{
        \input{figs/inter_list.tex}
    }
    \caption{$\mathcal{H}^2$~matrix, further decomposition of matrix $\widehat{M}_L$
    into sum of close and far interaction matrices $M_L$ and $F_{L-1}$. Close blocks
    are denoted by cyan color, far blocks are illustrated by blue color. White
    colored blocks are zero-valued.
    }
    \label{fig:inter_list}
\end{figure}
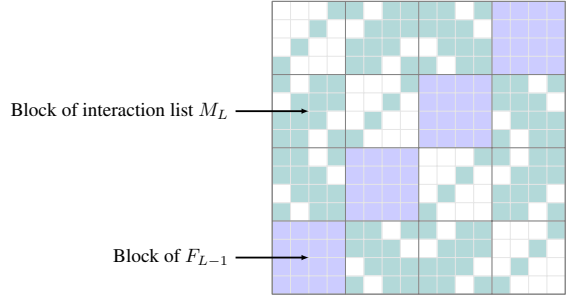

Take a look at receiver boxes that are close to source boxes on the current level.
We say that they set the {\bf interaction list}.
\begin{definition}
    A far block is in interaction list on level $l$ if the block that contains
    it on level $l-1$ is a close block, see Figure~\ref{fig:inter_list}.
    Blocks that correspond to interaction list on level $l$ constitute
    the interaction list matrix $M_{l+1}$.
\end{definition}

We separate $\widehat{M}_L$ into two parts: the interaction list matrix $M_L$,
and compressible matrix $F_{L-1}$, see Figure~\ref{fig:inter_list}.
\begin{equation}
\widehat{M}_L = M_L + F_{L-1}.
\end{equation}
Therefore,
\begin{equation}
F_L = U_L\widehat{M}_LV_L =  U_LM_LV_L + U_LF_{L-1}V_L
\end{equation}
Thus, we obtain:
\begin{equation}
w = Cq + U_LM_LV_Lq + U_LF_{L-1}V_Lq.
\end{equation}

\subsection{Compression at next levels}
\label{sec:lc}

For far blocks $F_{L-1}$ from level $l=L-1$, we compute the compression matrices
$U_{L-1}$ and $ V_{L-1}$, by the analogy to level $l=L$.
Repeating till the level $l=1$ we obtain the final formula:
\begin{equation}
    \label{eq:h2}
    w = Cq + \sum_{l=L}^{1} \left( U_L \ldots U_{L-l+1} \right) M_l
    \left( V_{L-l+1} \ldots V_{L} \right) q.
\end{equation}

\subsection{$\mathcal{H}^2$~matrix by vector product}
\label{sec:hmbvp}

Let us describe the $\mathcal{H}^2$~matrix by vector multiplication procedure.
Matrix $\mathcal{H}^2$~is given by equation~\eqref{eq:h2}

\begin{algorithm}[H]
\SetKwInOut{Input}{input}
\SetKwInOut{Output}{output}
\KwIn{Matrices $C$, $U_l$, $V_l$, $M_l$, $l=1 \ldots L$, vector $q$}
\KwOut{Vector $w$ defined by equation~\eqref{eq:h2}}
\SetAlgoLined
 $w = 0$\\
 $\hat{w}_l = 0, \; l=1\dots L$\\
 $w = w + Cq$\\
 \For{$l\gets1$ \KwTo $L$}{
    \For{$i\gets L$ \KwTo $L+1-l$}{
        $\hat{q}_i = V_i \hat{q}_{i+1}$ ($ \hat{q}_{L+1} = q$)
    }
    $\hat{w}_l = \hat{w}_l + M_l \hat{q}_l$\;
    \For{$i\gets l$ \KwTo $L$}{
        $\hat{w}_{i+1} = \hat{w}_{i+1} + U_i \hat{w}_{i}$ ($ \hat{w}_{L+1} = w$)
    }
}
 \caption{Pseudo-code for computation of $\mathcal{H}^2$~matrix by vector product}
 \label{alg:matvec}
\end{algorithm}

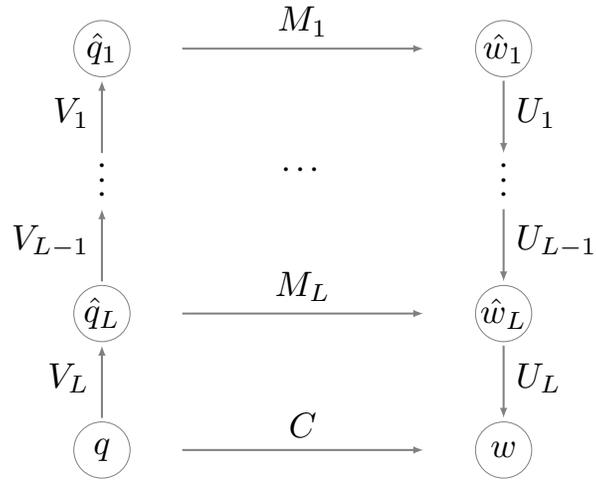
\begin{figure}[H]
\centering
    \resizebox{0.6\textwidth}{!}{
        \input{figs/matvec.tex}
    }
    \caption{Computational diagram of $\mathcal{H}^2$~matrix by vector product}
    \label{fig:h2_matvec}
\end{figure}

In Figure~\ref{fig:h2_matvec} the illustration of Algorithm~\ref{alg:matvec} is shown.
Vectors are placed in circles, arrows represent the matrix by vector product.

\begin{figure}[H]
\centering
    \resizebox{0.6\textwidth}{!}{
        \input{figs/nn.tex}
    }
\caption{Computational diagram of $\mathcal{H}^2$-based neural network}
\label{fig:hnn}
\end{figure}
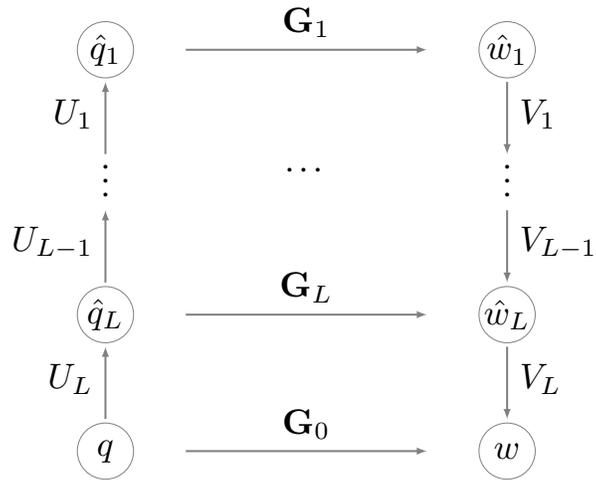

\end{document}

%% file: figs/close_mat.tex
\def\closeDots{{{ 0, 10,  5,  2,  1},
        {1, 11,  0,  3,  4},
        {4, 14,  1,  6,  5},
        {5, 15,  4,  7,  0},
        {2,  0,  7,  8,  3},
        {3,  1,  2,  9,  6},
        {6,  4,  3, 12,  7},
        {7,  5,  6, 13,  2},
        {8,  2, 13, 10,  9},
        {9,  3,  8, 11, 12},
       {12,  6,  9, 14, 13},
       {13,  7, 12, 15,  8},
       {10,  8, 15,  0, 11},
       {11,  9, 10,  1, 14},
       {14, 12, 11,  4, 15},
       {15, 13, 14,  5, 10}}}

\def\RandDots{{{0.48637352, 0.53586711, 0.87386436, 0.92523105, 0.19662421,
               0.11208547, 0.34445174, 0.14922475, 0.19785328, 0.63301906,
               0.25910754, 0.58501838, 0.88157351, 0.05393866, 0.31015342,
               0.08242764, 0.72402543, 0.61453043, 0.7415675 , 0.86744318},
              {0.94352635, 0.55941481, 0.66484183, 0.82699268, 0.16481487,
               0.52753435, 0.34134099, 0.56140683, 0.67282776, 0.93078796,
               0.30185688, 0.77235954, 0.7517697 , 0.93256297, 0.48900417,
               0.43782538, 0.12732212, 0.90118619, 0.56294379, 0.93998838}}}
\def\RandDotsT{{{0.8245779 , 0.76174813, 0.80744847, 0.97529873, 0.22327913,
               0.62019291, 0.15440187, 0.01456596, 0.82576358, 0.49669774,
               0.13688487, 0.55712594, 0.36549842, 0.84218367, 0.26413289,
               0.37919453, 0.72994362, 0.36359501, 0.63787649, 0.15465958},
              {0.56214302, 0.16206764, 0.26705701, 0.88304504, 0.566555  ,
               0.14174434, 0.61395399, 0.22522736, 0.56733802, 0.07072593,
               0.53008812, 0.93976755, 0.55426911, 0.05349953, 0.02133711,
               0.74527227, 0.94190781, 0.12909986, 0.84066342, 0.39037355}}}

%% file: figs/tree.tex
\tikz{
    \draw(0,-5.5)--(1,-4)--(1.8,-5.5)--(1,-4)--(2,-2.5)--(3,-4)--(2.2,-5.5)--(3,-4)--(4,-5.5);
    \path[fill=green!80!blue,draw=green] (0,-5.5) circle (1mm);
    \path[fill=green!80!blue,draw=green] (1,-4) circle (1mm);
    \path[fill=green!80!blue,draw=green] (1.8,-5.5) circle (1mm);
    \path[fill=green!80!blue,draw=green] (2,-2.5) circle (1mm);
    \path[fill=green!80!blue,draw=green] (3,-4) circle (1mm);
    \path[fill=green!80!blue,draw=green] (2.2,-5.5) circle (1mm);
    \path[fill=green!80!blue,draw=green] (4,-5.5) circle (1mm);
    \draw(7,-5.5)--(8,-4)--(8.8,-5.5)--(8,-4)--(9,-2.5)--(10,-4)--(9.2,-5.5)--(10,-4)--(11,-5.5);
    \path[fill=red!80!blue,draw=red] (7,-5.5) circle (1mm);
    \path[fill=red!80!blue,draw=red] (8,-4) circle (1mm);
    \path[fill=red!80!blue,draw=red] (8.8,-5.5) circle (1mm);
    \path[fill=red!80!blue,draw=red] (9,-2.5) circle (1mm);
    \path[fill=red!80!blue,draw=red] (10,-4) circle (1mm);
    \path[fill=red!80!blue,draw=red] (9.2,-5.5) circle (1mm);
    \path[fill=red!80!blue,draw=red] (11,-5.5) circle (1mm);
    \draw(0,-7)--(3,-7)--(3,-7.2)--(0,-7.2) -- (0,-7);
    \draw (1.5,-8) node[above]{\huge{$x$}} ;
    \draw [<-,  red] (2,-5.9) --(1.5,-6.6)   ;
    \draw (1.2,-6.4) node [above]  {};
    \draw (7.5,-7)--(10.5,-7)--(10.5,-7.2)--(7.5,-7.2) -- (7.5,-7);
    \draw (9,-8) node [above]  {\huge{ $y$}} ;
    \draw [<-,  red] (8.5,-6.6) -- (9,-5.9);
    \draw (8.3,-6.4) node [above]  {}; 
    \draw (11.9,-4.5) node [above]  {\Large{$U$}} ;
    \draw [<-,  red] (-1.5,-2.7) -- (-1.5,-5.7);
    \draw (5.5,-4.4) node [above]  {\Large{$M$}} ;
    \draw [<-,  red] (6.3,-4.4) -- (4.6,-4.4);
    \draw (-1.1,-4.5) node [above]  {\Large{$V$}} ;
   \draw [<-, red] (12.3,-5.7) -- (12.3,-2.7);
   \draw (5,-7) node [above]  {\Large{$C$}} ;
    \draw [<-,  red] (6,-7) -- (4,-7);
    \draw (2,-4.7) node [above]  {\huge{$\hat{x}$}} ;
    \draw (9,-4.7) node [above]  {\huge{$\hat{y}$}} ;
}

%% file: figs/nn.tex
\pgfmathsetmacro{\r}{0.7}
\pgfmathsetmacro{\scale}{2.5}
\pgfmathsetmacro{\ps}{10}
    \tikz{
        \draw[gray] (0*\ps,0*\ps) circle (\r cm);
        \node [scale=\scale] at (0*\ps,0*\ps) {$q$};
        \draw[draw=gray] (1*\ps,0*\ps) circle (\r cm);
        \node [scale=\scale] at (1*\ps,0*\ps) {$w$};
        \draw[-latex,line width=0.5mm ,draw=gray] (0.2*\ps ,0*\ps) -- (0.8*\ps,0*\ps);
        \node [scale=\scale, above] at (0.5*\ps,0*\ps) {${\bf G}_0$};
        \draw[-latex,line width=0.5mm, draw=gray] (0*\ps,0.08*\ps) -- (0*\ps,0.26*\ps);
        \node [scale=\scale, left] at (0*\ps,0.18*\ps) {$U_L$};
        \draw[latex-,line width=0.5mm, draw=gray] (1*\ps,0.08*\ps) -- (1*\ps,0.26*\ps);
        \node [scale=\scale, right] at (1*\ps,0.18*\ps) {$V_L$};

        \begin{scope}[shift={(0,3.4)}]
            \draw[gray] (0*\ps,0*\ps) circle (\r cm);
            \node [scale=\scale] at (0*\ps,0*\ps) {$\hat{q}_L$};
            \draw[draw=gray] (1*\ps,0*\ps) circle (\r cm);
            \node [scale=\scale] at (1*\ps,0*\ps) {$\hat{w}_L$};
            \draw[-latex,line width=0.5mm, draw=gray] (0.2 *\ps,0*\ps) -- (0.8*\ps,0*\ps);
            \node [scale=\scale, above] at (0.5*\ps,0*\ps) {${\bf G}_L$};
            \draw[-latex,line width=0.5mm, draw=gray] (0*\ps,0.08*\ps) -- (0*\ps,0.26*\ps);
            \node [scale=\scale, left] at (0*\ps,0.18*\ps) {$U_{L-1}$};
            \draw[latex-,line width=0.5mm, draw=gray] (1*\ps,0.08*\ps) -- (1*\ps,0.26*\ps);
            \node [scale=\scale, right] at (1*\ps,0.18*\ps) {$V_{L-1}$};
        \end{scope}
        \begin{scope}[shift={(0,6.6)}]
            \draw[-latex,line width=0.5mm, draw=gray] (0*\ps,0.08*\ps) -- (0*\ps,0.26*\ps);
            \node [scale=\scale, left] at (0*\ps,0.18*\ps) {$U_{1}$};
            \draw[latex-,line width=0.5mm, draw=gray] (1*\ps,0.08*\ps) -- (1*\ps,0.26*\ps);
            \node [scale=\scale, right] at (1*\ps,0.18*\ps) {$V_{1}$};
        \end{scope}
        \begin{scope}[shift={(0,1*\ps)}]
            \draw[gray] (0*\ps,0*\ps) circle (\r cm);
            \node [scale=\scale] at (0*\ps,0*\ps) {$\hat{q}_1$};
            \draw[draw=gray] (1*\ps,0*\ps) circle (\r cm);
            \node [scale=\scale] at (1*\ps,0*\ps) {$\hat{w}_1$};
            \draw[-latex,line width=0.5mm,draw=gray ] (0.2*\ps ,0*\ps) -- (0.8*\ps,0*\ps);
            \node [scale=\scale, above] at (0.5*\ps,0*\ps) {${\bf G}_1$};
        \end{scope}

        \node [scale=\scale] at (0.5*\ps,0.7*\ps) {\dots};
        \node [scale=\scale] at (0*\ps,0.7*\ps) {\vdots};
        \node [scale=\scale] at (1*\ps,0.7*\ps) {\vdots};
        \node [scale=\scale] at (5,-1.5) {~};
    }

%% file: figs/rte_table.tex
 \begin{tabular}{||c c c c c||}
 \hline
 $N$ & h2nn &  MNN-$\mathcal{H}^2$-Mix &  MNN-$\mathcal{H}^2$-LC & conv\\
 \hline\hline
 320 &0.02350/0.02346 &   0.01197/0.01196&   0.00366/0.00383& 0.19538/0.19359\\
 640 &  0.02676/0.02687 &  0.06179/0.06381&   0.06271/0.06692&0.19683/0.19671\\
 1280 &   0.03240/0.03233 &  0.12397/0.125108&   0.16229/0.16438& 0.19686/0.19552\\ [1ex]
 \hline
\end{tabular}

%% file: figs/ce_table.tex
 \begin{tabular}{||c c c c c||}
 \hline
 $N$ & h2nn &  MNN-$\mathcal{H}^2$-Mix &  MNN-$\mathcal{H}^2$-LC& conv\\
 \hline\hline
 320 &  0.00790/0.00789 &   0.00857/0.00858 &   0.00218/0.00218&0.26599/0.26600\\
 640 &  0.00424/0.00425 &  0.14413/0.14592 &   0.10721/0.11022&0.27551/0.27860\\
 1280 &   0.01464/0.00425 &  0.19069/0.19070 &   0.22265/0.22287&0.33206/0.33301 \\ [1ex]
 \hline
\end{tabular}

%% file: figs/separation_property.tex
\begin{tikzpicture}[scale=0.3]
    \pgfmathsetmacro{\r}{0.1}
    \pgfmathsetmacro{\n}{19}
    \foreach \i in {0,...,\n}{
            \draw[draw=gray,fill=gray] (\RandDots[0][\i]*5,\RandDots[1][\i]*5) circle (0.1 cm);
    }
    \draw [gray, opacity=1.0] (0,0) rectangle ++(5,5);
    \node [scale=1] at (2.5,-1) {$x$};
    \begin{scope}[shift={(15,0)}]
        \foreach \i in {0,...,\n}{
                \draw[draw=green,fill=green] (\RandDotsT[1][\i]*5,\RandDots[0][\i]*5) circle (0.1 cm);
        }
        \draw [green, opacity=1.0] (0,0) rectangle ++(5,5);
        \node [scale=1] at (2.5,-1) {$y$};
    \end{scope}

    \draw[latex-latex,line width=0.3mm,draw=gray,fill=gray] (6,2.5) -- (14,2.5);
\end{tikzpicture}

%% file: figs/grid_a.tex
\begin{tikzpicture}[scale=\size]
        \pgfmathsetmacro{\r}{0.1}
        \pgfmathsetmacro{\n}{8}
        \foreach \j in {1,...,\n}{
            \foreach \i in {1,...,\n}{
                \draw[draw = gray,fill=gray] (\i,\j+\r) arc (90:270:\r cm);
                \draw[draw =green,fill=green] (\i,\j-\r) arc (-90:90:\r cm);
            }
        }
        \draw[draw=blue!50!green] (0.5,0.5) rectangle (\n+0.5,\n+0.5);
        \draw[draw=blue!50!green] (0.5,\n/2+0.5) -- (\n+0.5,\n/2+0.5);
        \draw[draw=blue!50!green] (\n/2+0.5,0.5) -- (\n/2+0.5,\n+0.5);
        \draw[draw=blue!50!green] (0.5,0.5+3*\n/4) -- (0.5+\n/2,0.5+3*\n/4);
        \draw[draw=blue!50!green] (\n/4+0.5,0.5+\n/2) -- (\n/4+0.5,0.5+\n);

\end{tikzpicture}

%% file: figs/grid_b.tex
\pgfmathsetmacro{\r}{0.1}
\pgfmathsetmacro{\n}{8}
\begin{tikzpicture}[scale=\size]
    \draw [blue!50!green,fill=blue!50!green,opacity=0.2](3.5,5.5) circle (2.3cm);
    \draw[draw=blue!50!green, fill=blue!50!green, opacity=0.3] (0.5,4.5) rectangle (2.5,6.5);
    \draw [draw=blue!50!green, fill = blue!50!green] (1.5,5.5) circle (\r cm/2);
    \draw[draw=blue!50!green, fill=blue!50!green, opacity=0.3] (2.5,4.5) rectangle (4.5,6.5);
    \draw [draw=blue!50!green, fill = blue!50!green] (3.5,5.5) circle (\r cm/2);
    \draw[draw=blue!50!green, fill=blue!50!green, opacity=0.3] (4.5,4.5) rectangle (6.5,6.5);
    \draw [draw=blue!50!green, fill = blue!50!green] (5.5,5.5) circle (\r cm/2);
    \draw[draw=blue!50!green, fill=blue!50!green, opacity=0.3] (2.5,2.5) rectangle (4.5,4.5);
    \draw [draw=blue!50!green, fill = blue!50!green] (3.5,3.5) circle (\r cm/2);
    \draw[draw=blue!50!green, fill=blue!50!green, opacity=0.3] (2.5,6.5) rectangle (4.5,8.5);
    \draw [draw=blue!50!green, fill = blue!50!green] (3.5,7.5) circle (\r cm/2);
        \foreach \j in {1,...,\n}{
            \foreach \i in {1,...,\n}{
                \draw[draw = gray,fill=gray] (\i,\j+\r) arc (90:270:\r cm);
                \draw[draw = green,fill=green] (\i,\j-\r) arc (-90:90:\r cm);
            }
        }
        \draw[step=2.0,blue!50!green,xshift=0.5cm,yshift=0.5cm] (0,0) grid (\n,\n);
\end{tikzpicture}

%% file: figs/mat_A.tex
\tikz{
     \fill[blue!20!white]  (0,0) rectangle (16,16);
    \foreach \i in {0,...,15}
    {
        \foreach \j in {0,...,4}
        {
            \pgfmathtruncatemacro{\r}{int(\closeDots[\i][\j])}
            \pgfmathtruncatemacro{\ii}{int(\closeDots[\i][0])}
            \draw [blue!50!white,fill=blue!50!white,opacity=1.0](\r,15-\ii) rectangle ++(1,1);
        }
    }
    \draw[step=1cm,gray!20!white] (0,0) grid (16,16);
    \node [scale=3,left] at (-2,14.5) {Block of close matrix $C$};
    \draw[-latex,line width=1mm ] (-2,14.5) -- (0.5,14.5);
    \node [scale=3,left] at (-2,11.5) {Block of far matrix $F_L$};
    \draw[-latex,line width=1mm ] (-2,11.5) -- (0.5,11.5);
    \draw[white] (30,0) rectangle ++(1,1);
}

%% file: figs/bottom_level_compression.tex
\tikz{
     \fill[blue!20!white]  (0,0) rectangle (16,16);
    \foreach \i in {0,...,15}
    {
        \foreach \j in {0,...,4}
        {
            \pgfmathtruncatemacro{\r}{int(\closeDots[\i][\j])}
            \pgfmathtruncatemacro{\ii}{int(\closeDots[\i][0])}
            \draw [white,fill=white,opacity=1.0](\r,15-\ii) rectangle ++(1,1);
        }
    }

    \draw[white,fill=white] (0,15) rectangle (16,15.5);
    \draw[white,fill=white] (3.5,0) rectangle (4,16);
    \draw[step=1cm,gray!20!white] (0,0) grid (16,16);

    \draw[green!50!white,fill=green!50!white](-2.5, 15) rectangle ++(0.5,1);
    \draw[green!50!white,fill=green!50!white](3,17) rectangle ++(1,0.5);

    \node [scale=2, above] at (3.5,17.5) {$\widetilde{V}_j$};
    \node [scale=2, left] at (-2.5,15.5) {$\widetilde{U}_i$};
    \node [scale=2] at (1.5,16.5) {$\widetilde{F}_{ij}$};
    \draw[-latex,line width=1mm ] (2,16.5) -- (3.25,15.75);
}

%% file: figs/L_lvl_compr.tex
\resizebox{\textwidth}{!}{
\tikz{
    \foreach \i in {0,...,15}
    {
        \foreach \j in {0,...,4}
        {
            \pgfmathtruncatemacro{\r}{int(\closeDots[\i][\j])}
            \pgfmathtruncatemacro{\ii}{int(\closeDots[\i][0])}
            \draw [blue!50!white,fill=blue!50!white,opacity=1.0](\r,15-\ii) rectangle ++(1,1);
        }
    }
    \draw[step=1cm,gray!20!white] (0,0) grid (16,16);
    \draw[gray!80!white] (0,0) rectangle (16,16);
    \node [scale=4] at (8,-2) { $C$};

    \begin{scope}[shift={(20,0)}]

            \foreach \i in {0,...,16}
            {
                \draw[gray!20!white] (0,\i)--(8,\i);
                \draw[gray!20!white] (\i*0.5,0) -- (\i*0.5,16);
            }
            \foreach \i in {0,...,15}
            {
            \draw[green!50!white,fill=green!50!white](\i*0.5,15-\i) rectangle ++(0.5,1);
            }
            \node [scale=4] at (4,-2) {$U_L$ };
            \draw[gray!80!white] (0,0) rectangle (8,16);
    \end{scope}
    \begin{scope}[shift={(30,0)}]
            \fill[blue!20!white]  (0,8) rectangle (8,16);
            \foreach \i in {0,...,15}
            {
                \foreach \j in {0,...,4}
                {
                    \pgfmathtruncatemacro{\r}{int(\closeDots[\i][\j])}
                    \pgfmathtruncatemacro{\ii}{int(\closeDots[\i][0])}
                    \draw [white,fill=white,opacity=1.0](\r*0.5,15.5-\ii*0.5) rectangle ++(0.5,0.5);
                }
            }
            \draw[step=1cm,gray!20!white] (0,8) grid (8,16);
            \draw[gray!80!white] (0,8) rectangle (8,16);
            \node [scale=4] at (4,6) {  $\widehat{M}_L$ };
    \end{scope}
    \begin{scope}[shift={(40,0)}]
        \foreach \i in {0,...,16}
        {
            \draw[gray!20!white] (0,\i*0.5+8)--(16,\i*0.5+8);
            \draw[gray!20!white] (\i,8) -- (\i,16);
        }
        \foreach \i in {0,...,15}
        {
        \draw[green!50!white,fill=green!50!white](\i,15.5-\i*0.5) rectangle ++(1,0.5);
        }
        \node [scale=4] at (8,6) {  $V_L$ };
        \draw[gray!80!white] (0,8) rectangle (16,16);
    \end{scope}
    \node [scale=4] at (18,8) {  + };
    \node [scale=4] at (29,12) {  * };
    \node [scale=4] at (39,12) {  * };
}
}

%% file: figs/inter_list.tex
    \tikz{
        \fill[blue!50!green, opacity=0.3]  (0,0) rectangle (16,16);
        \foreach \i in {0,...,3}{
            \fill[blue!20!white]  (\i*4,\i*4) rectangle ++(4cm,4cm);
        }
        \foreach \i in {0,...,15}
        {
            \foreach \j in {0,...,4}
            {
                \pgfmathtruncatemacro{\r}{int(\closeDots[\i][\j])}
                \pgfmathtruncatemacro{\ii}{int(\closeDots[\i][0])}
                \draw [white,fill=white,opacity=1.0](\r,15-\ii) rectangle ++(1,1);
            }
        }
        \draw[step=1cm,gray!20!white] (0,0) grid (16,16);
        \draw[gray!80!white] (0,0) rectangle (16,16);
        \draw[step=4cm,gray,line width=0.34mm ] (0,0) grid (16,16);
        \node [scale=3,left] at (-2,2) {Block of $F_{L-1}$};
        \draw[-latex,line width=1mm ] (-2,2) -- (2,2);
        \node [scale=3,left] at (-2,10) {Block of interaction list $M_L$};
        \draw[-latex,line width=1mm ] (-2,10) -- (2,10);
        \draw[white] (30,0) rectangle ++(1,1);
    }

%% file: figs/matvec.tex
\pgfmathsetmacro{\r}{0.7}
\pgfmathsetmacro{\scale}{2.5}
\pgfmathsetmacro{\ps}{10}
    \tikz{
        \draw[gray] (0*\ps,0*\ps) circle (\r cm);
        \node [scale=\scale] at (0*\ps,0*\ps) {$q$};
        \draw[draw=gray] (1*\ps,0*\ps) circle (\r cm);
        \node [scale=\scale] at (1*\ps,0*\ps) {$w$};
        \draw[-latex,line width=0.5mm ,draw=gray] (0.2*\ps ,0*\ps) -- (0.8*\ps,0*\ps);
        \node [scale=\scale, above] at (0.5*\ps,0*\ps) {$C$};
        \draw[-latex,line width=0.5mm, draw=gray] (0*\ps,0.08*\ps) -- (0*\ps,0.26*\ps);
        \node [scale=\scale, left] at (0*\ps,0.18*\ps) {$V_L$};
        \draw[latex-,line width=0.5mm, draw=gray] (1*\ps,0.08*\ps) -- (1*\ps,0.26*\ps);
        \node [scale=\scale, right] at (1*\ps,0.18*\ps) {$U_L$};

        \begin{scope}[shift={(0,3.4)}]
            \draw[gray] (0*\ps,0*\ps) circle (\r cm);
            \node [scale=\scale] at (0*\ps,0*\ps) {$\hat{q}_L$};
            \draw[draw=gray] (1*\ps,0*\ps) circle (\r cm);
            \node [scale=\scale] at (1*\ps,0*\ps) {$\hat{w}_L$};
            \draw[-latex,line width=0.5mm, draw=gray] (0.2 *\ps,0*\ps) -- (0.8*\ps,0*\ps);
            \node [scale=\scale, above] at (0.5*\ps,0*\ps) {$M_L$};
            \draw[-latex,line width=0.5mm, draw=gray] (0*\ps,0.08*\ps) -- (0*\ps,0.26*\ps);
            \node [scale=\scale, left] at (0*\ps,0.18*\ps) {$V_{L-1}$};
            \draw[latex-,line width=0.5mm, draw=gray] (1*\ps,0.08*\ps) -- (1*\ps,0.26*\ps);
            \node [scale=\scale, right] at (1*\ps,0.18*\ps) {$U_{L-1}$};
        \end{scope}
        \begin{scope}[shift={(0,6.6)}]
            \draw[-latex,line width=0.5mm, draw=gray] (0*\ps,0.08*\ps) -- (0*\ps,0.26*\ps);
            \node [scale=\scale, left] at (0*\ps,0.18*\ps) {$V_{1}$};
            \draw[latex-,line width=0.5mm, draw=gray] (1*\ps,0.08*\ps) -- (1*\ps,0.26*\ps);
            \node [scale=\scale, right] at (1*\ps,0.18*\ps) {$U_{1}$};
        \end{scope}
        \begin{scope}[shift={(0,1*\ps)}]
            \draw[gray] (0*\ps,0*\ps) circle (\r cm);
            \node [scale=\scale] at (0*\ps,0*\ps) {$\hat{q}_1$};
            \draw[draw=gray] (1*\ps,0*\ps) circle (\r cm);
            \node [scale=\scale] at (1*\ps,0*\ps) {$\hat{w}_1$};
            \draw[-latex,line width=0.5mm,draw=gray ] (0.2*\ps ,0*\ps) -- (0.8*\ps,0*\ps);
            \node [scale=\scale, above] at (0.5*\ps,0*\ps) {$M_1$};
        \end{scope}

        \node [scale=\scale] at (0.5*\ps,0.7*\ps) {\dots};
        \node [scale=\scale] at (0*\ps,0.7*\ps) {\vdots};
        \node [scale=\scale] at (1*\ps,0.7*\ps) {\vdots};
    }